

\input epsf.tex

\def\2{{1\over 2}}

\def\d{\delta}
\def\a{\alpha}
\def\b{\beta}
\def\g{\gamma}

\def\e{\epsilon}
\def\l{\lambda}

\def\fun#1#2#3{#1\colon #2\rightarrow #3}

\def\frac#1#2{{{#1} \over {#2}}}

\def\sqr{\sqrt}
\def\st{\;\colon\;}
\def\tends{\rightarrow}
\def\weak{\rightharpoonup}
\def\sc#1#2{(#1\vert #2)}

\def\dr{ {\rm d} }

\def\C{{\bf C}}

\def\R{{\bf R}}

\def\thm#1{\vskip 1 pc\noindent{\bf Theorem #1.\quad}\sl}
\def\lem#1{\vskip 1 pc\noindent{\bf Lemma #1.\quad}\sl}
\def\prop#1{\vskip 1 pc\noindent{\bf Proposition #1.\quad}\sl}

\def\proof{\rm\vskip 1 pc\noindent{\bf Proof.\quad}}
\def\fin{\par\hfill $\backslash\backslash\backslash$\vskip 1 pc}
\def\txt#1{\quad\hbox{#1}\quad}

\def\L{{\cal L}}

\def\2{\frac{1}{2}}

\def\part{{\partial_{x}}}

\def\sc{{\cal S}}
\def\kc{{\cal K}}
\def\bc{{\cal B}}

\def\ec{{\cal E}}

\def\ac{{\cal A}}

\def\uc{{\cal U}}
\def\dc{{\cal D}}
\def\oc{{\cal O}}

\def\kc{{\cal K}}
\def\uc{{\cal U}}

\def\rc{{\cal R}}



\baselineskip= 17.2pt plus 0.6pt
\font\titlefont=cmr17
\centerline{\titlefont Harmonic immersions of the Sierpinski gasket}
\vskip 1 pc
\centerline{\titlefont into the hyperbolic plane.}
\vskip 4pc
\font\titlefont=cmr12
\centerline{         \titlefont {Ugo Bessi}\footnote*{{\rm 
Dipartimento di Matematica, Universit\`a\ Roma Tre, Largo S. 
Leonardo Murialdo, 00146 Roma, Italy.}}   }{}\footnote{}{
{{\tt email:} {\tt bessi@matrm3.mat.uniroma3.it} Work partially supported by the PRIN2023 grant "Stability in Hamiltonian Systems and beyond"}} 
\vskip 0.5 pc
 
\par
\vskip 2pc
\centerline{\bf Abstract} 

Many fractals $G$ admit a harmonic immersion into $\R^n$, i.e. an immersion which minimises a natural energy under fixed boundary conditions; we look for harmonic immersions of the Sierpinski gasket into the hyperbolic plane. We show that, given any three points $\tilde A$, $\tilde B$, $\tilde C$ in the hyperbolic plane there is a harmonic map bringing the three points $A$, $B$, $C$ of the boundary of the gasket to $\tilde A$, $\tilde B$, $\tilde C$ respectively. Moreover, if the points 
$\tilde A$, $\tilde B$, $\tilde C$ are sufficiently close in the hyperbolic distance, then the harmonic map is unique and depends differentiably on $\tilde A$, $\tilde B$, $\tilde C$. Lastly, we show that, if the harmonic map $\phi$ is injective, then it brings geodesics of the gasket $G$ into geodesics of 
$\phi(G)$.

\vskip 2 pc
\centerline{\bf  Introduction}
\vskip 1 pc

Many fractals $G\subset\R^n$ are defined in the following way ([15], [21]): we are given finitely many invertible, affine contractions 
$$\fun{F_i}{\R^n}{\R^n},\qquad i\in\{ 1,\dots,p \}$$
and $G$ is the unique non-empty, compact set of $\R^n$ which satisfies 
$$G=\bigcup_{i=1}^p F_i(G)  .  \eqno (1)$$
For several of these fractals, one can define a bilinear form on $G$ which is formally similar to Dirichlet's energy on $\R^n$. More precisely, it is possible to define, in a natural way, a Borel probability measure $\kappa$ on $G$ (a.k.a. Kusuoka's measure, [18]) and a Borel field $E_x$ of symmetric, semidefinite positive matrices from the dual space of 
$\R^n$ to itself which induce the bilinear form 
$$\fun{\ec^\dr}{C^1(\R^n,\R)\times C^1(\R^n,\R)}{\R}$$
$$\ec^\dr( u, v)\colon=
\int_G(\dr u(x)E_x,\dr v(x))\dr\kappa(x)  .  \eqno (2)$$
The inner product in the formula above is the natural one in the dual of $\R^n$ and $\dr$ is the standard differential. 

Following [14] and [16], we consider $E_x$ as a Riemannian tensor on the "cotangent space" of $G$ and $\kappa$ as a Riemannian volume on $G$. We have written 
$\dr u(x)E_x$ because we adopt the convention that row vectors are in the cotangent space and column vectors are in the tangent one. 

Having $\ec^\dr$, we say that a function $u$ is harmonic if 
$$\ec^\dr( u,\phi)=0$$
for all $\phi\in C^1(\R^n,\R)$ which vanish on the "boundary" of the fractal. 

It is standard (see for instance [8]) that, when the coordinate functions are harmonic, the form $\ec^\dr$ is closable in $L^2(G,\kappa)$ and thus it extends to a Dirichlet form defined on $\dc(\ec^\dr)\subset L^2(G,\kappa)$ with $C^1(\R^n,\R)\subset \dc(\ec^\dr)$. Some famous cases of this situation are the harmonic Sierpinski gasket ([14], [16]), the 
level-$n$ Sierpinski gasket ([10], [25]) and the stretched Sierpinski gasket ([2], [7]). 

Saying that the coordinate functions are harmonic is equivalent to say that, under suitable boundary conditions, the identity map minimises the energy of the immersion into $\R^n$. The energy $\tilde\ec$ of the immersion $\fun{\phi}{G\subset\R^n}{\R^n}$ is defined in the standard way: if 
$$\phi(x_1,\dots,x_n)=(
u_1(x_1,\dots,x_n),\dots,u_n(x_1,\dots,x_n)
)$$
then, for the form $\ec^\dr$ of (2), 
$$\tilde\ec(\phi)\colon=
\ec^\dr( u_1,u_1)+\dots +\ec^\dr( u_n, u_n)  .  $$
Harmonic immersions of a Riemannian manifold into another one have been studied extensively since [11]; a natural question is whether there are energy-minimising immersions of fractals into manifolds different from $\R^n$. 

In this paper we study a particular case and embed the harmonic Sierpinski gasket $G$ into the hyperbolic plane 
$$H_+=\{
x+iy\st y>0
\}   \eqno (3)$$
endowed with the Poincar\'e\ metric. We briefly recall that, if $v$ and $w$ are in the tangent space $T_{(x,y)}H_+$, then their hyperbolic inner product is given by 
$$(v,w)_{Hyp}=\frac{(v,w)}{y^2}  $$
where $(\cdot,\cdot)$ denotes the Euclidean inner product. 

Now we can follow [11] and define the energy of an immersion 
$$\fun{\phi=u+iv}{G}{H_+}  \eqno (4)$$
with $u,v\in\dc(\ec)$ as 
$$\tilde\ec(\phi)=
\int_G\frac{1}{v^2}\left[
\left( \dr u E_x,\dr u \right)+\left( \dr v E_x,\dr v \right)
\right]\dr\kappa  .  \eqno (5)$$
We recall from [14] that the function in $\dc(\ec^\dr)$ are continuous, which allows us to define the boundary conditions in the following way. 

\vskip 1pc

\noindent{\bf Definitions.} Let $\tilde A,\tilde B,\tilde C$ be three assigned points in $H_+$; let 
$\{ A,B,C \}\subset\R^2$ be the vertices of the equilateral triangle of figure 1 below; they are the "boundary" of the harmonic gasket. We define $\ac$ as the set of the maps 
$\phi=u+iv$ as in (4) such that 

\noindent 1) $u,v\in\dc(\ec^\dr)$  .  

\noindent 2) $\phi(A)=\tilde A$, $\phi(B)=\tilde B$, $\phi(C)=\tilde C$. 

For the energy $\tilde\ec$ of (5) we define 
$$\a=\inf_{\phi\in\ac}\tilde\ec(\phi)  .  \eqno (6)$$
We say that the map $\phi\in\ac$ is harmonic if $\tilde\ec(\phi)=\a$. 

\vskip 1pc

We shall prove the following theorem. 

\thm{1} 1) Given any three points $\tilde A,\tilde B,\tilde C\in H_+$, the number $\a$ of (6) is finite and there is at least one harmonic map, i.e. a map 
$\phi\in\ac$ with $\tilde\ec(\phi)=\a$.  

\noindent 2) There is $r_0>0$ such that, if $\tilde A$, $\tilde B$, $\tilde C$ are contained in a ball of hyperbolic radius smaller that $r_0$, then the minimal map $\phi$ is unique and depends differentiably on $\tilde A$, $\tilde B$, $\tilde C$. 

\noindent 3) If the immersion $\phi$ is injective, then it is totally geodesic: in other words, curves in $G$ which locally are geodesics are brought into curves in $\phi(G)$ which locally are geodesics.

\rm

\vskip 1pc

A few comments. First, when the target space is $\R^n$, the maps $F_i$ of (1) are closely related to the harmonic extension maps on the pre-fractals ([21], [23]); they are affine because the harmonic extension is linear. This is no more the case in our setting, which is non-linear, but we can modify (1) in the following way. We call 
$\fun{\phi(\tilde A,\tilde B,\tilde C)}{G}{H_+}$ the harmonic map with boundary conditions 
$\tilde A$, $\tilde B$, $\tilde C$; then, for the points $a$, $b$, $c$ of figure 1 below, 
$$\phi(\tilde A,\tilde B,\tilde C)(G)=
\phi(\tilde A,\phi(\tilde A,\tilde B,\tilde C)(c),\phi(\tilde A,\tilde B,\tilde C)(b))(G)\cup$$
$$\phi(\phi(\tilde A,\tilde B,\tilde C)(c),\tilde B,\phi(\tilde A,\tilde B,\tilde C)(a))(G)\cup
\phi(\phi(\tilde A,\tilde B,\tilde C)(b),\phi(\tilde A,\tilde B,\tilde C)(a),\tilde C)(G)  .  $$ 
We won't give a proof of this fact, which says that the image of $G$ under $\phi$ is the union of the images of the three harmonic maps which bring the triple $(A,B,C)$ into the triples 
$$(\tilde A,\phi(\tilde A,\tilde B,\tilde C)(c),\phi(\tilde A,\tilde B,\tilde C)(b)),
(\phi(\tilde A,\tilde B,\tilde C)(c),\tilde B,\phi(\tilde A,\tilde B,\tilde C)(a)),
(\phi(\tilde A,\tilde B,\tilde C)(b),\phi(\tilde A,\tilde B,\tilde C)(a),\tilde C)  .  $$

We also note that we do not know whether a harmonic map $\phi$ is injective, even if the points 
$\tilde A,\tilde B,\tilde C$ do not lie on the same geodesic of $H_+$, and that's the reason for the additional hypothesis in point 3) of theorem 1. When instead of a fractal $G$ we have a manifold a stronger fact holds, namely that being harmonic is equivalent to being totally geodesic. 

Since the energy $\tilde\ec$ is non linear, the compactness needed for the existence of the minimum will follow from the Ascoli-Arzel\`a\ theorem; to apply it we shall need some delicate results of [1], [5] and [17], which imply that the domain of $\ec^\dr$ embeds in a space of H\"older functions.

The paper is organised as follows: in section 1 we recall the construction of the harmonic gasket, in section 2 we recall some properties of the energy $\ec$ of (2). In section 3 we recall from [1], [5], and [17] some facts about function spaces defined on $G$, in section 4 we prove point 1) of theorem 1, in section 5 we prove point 2) and in section 6 we prove point 3).

\vskip 2pc

\centerline{\bf \S 1}
\centerline{\bf The harmonic gasket and its energy}
\vskip 1pc

\noindent{\bf The harmonic Sierpinski gasket.} We briefly recall the definition and some properties of the Harmonic Sierpinski Gasket in $\R^2$; we refer the reader to [14] and [16] for first hand information. We set
$$T_1=\left(
\matrix{
\frac{3}{5},&0\cr
0,&\frac{1}{5}
}
\right)  ,  \quad
T_2=\left(
\matrix{
\frac{3}{10},&\frac{\sqr 3}{10}\cr
\frac{\sqr 3}{10},&\frac{1}{2}
}
\right)  ,\quad
T_3=\left(
\matrix{
\frac{3}{10},&-\frac{\sqr 3}{10}\cr
-\frac{\sqr 3}{10},&\frac{1}{2}
}
\right)   ,  $$
$$A=\left(
\matrix{
0\cr
0
}
\right),\qquad
B=\left(
\matrix{
1\cr
\frac{1}{\sqrt 3}
}
\right) ,\qquad
C=\left(
\matrix{
1\cr
-\frac{1}{\sqrt 3}
}
\right)  $$
and 
$$F_1(x)=T_1(x),\quad F_2(x)=
B   +T_2\left(
x-B
\right),\quad
F_3(x)=
C   +T_3\left(
x-C
\right)  .  $$
Referring to figure 1 below, $F_1$ brings the equilateral triangle $ABC$ into $Abc$; 
$F_2$ brings $ABC$ into $Bac$ and $F_3$ brings $ABC$ into $Cba$.

\input miniltx %
\input graphicx.sty %
\includegraphics[scale=0.5]{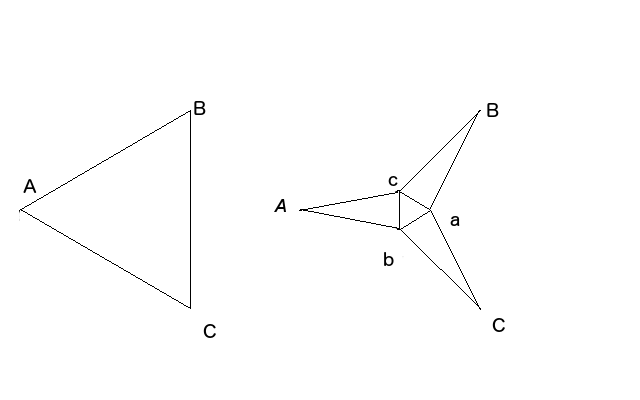}

\centerline{Figure 1}

\vskip 1pc

\noindent{\bf The bilinear form.} We recall a standard fact from [14] (for alternative approaches using dynamical systems, see [12], [19], [20] or [6].) Let $\Lambda^1(\R^2)$ denote the dual space of $\R^2$; we endow it with the inner product induced by the standard inner product of $\R^2$, i.e. 
$$(\dr x,\dr x)=1=(\dr y,\dr y)
\txt{and}
(\dr x,\dr y)=0  .  $$
In coordinates, for us the vectors of $\R^2$ will be column vectors and those of 
$\Lambda^1(\R^2)$ will be row vectors. 

We denote by $M$ the space of symmetric matrices from $\Lambda^1(\R^2)$ to itself. Then we can define, in a natural way, a Borel probability measure $\kappa$ on $G$ (which is called Kusuoka's measure, [18]) and a Borel function $\fun{E}{G}{M}$; it is known ([14]) that $E_x$ is a one-dimensional projection for $\kappa$-a.e. $x\in G$. 

These objects induce a bilinear form 
$$\fun{\ec}{C(G,\Lambda^1(\R^2))\times C(G,\Lambda^1(\R^2))}{\R}$$
$$\ec(u,v)=\int_G(u(x)E_x,v(x))\dr\kappa(x)  .  \eqno (1.1)$$
Note again that the symmetric matrix $E_x$ appears on the right because $u(x)$ is a row vector. 

It turns out ([14], but see also the next section) that the bilinear form of (2), i.e. 
$$\fun{\ec^\dr}{C^1(\R^2,\R)\times C^1(\R^2,\R)}{\R}$$
$$\ec^\dr(f,g)\colon=\ec(\dr f,\dr g)   \eqno (1.2)$$
is closable in $L^2(G,\kappa)$; its closure is a local Dirichlet form on $G$. 

\noindent{\bf The Lebesgue space of one-forms.} We call $\L$ be the space of the Borel functions $\fun{u}{G}{\Lambda^1(\R^2)}$ such that $u(x)\in {\rm Rank}(E_x)$ for $\kappa$-a.e. $x\in G$ and we define
$$L^2(G,\Lambda^1(G),\kappa)\colon=\{
u\in\L\st\ec(u,u)<+\infty
\}  .  $$
We forego the easy proof that $L^2(G,\Lambda^1(G),\kappa)$ is a Hilbert space for the natural inner product 
$$(u,v)_{L^2}\colon=\ec(u,v)=
\int_G(u(x)E_x,v(x))\dr\kappa(x)  .    \eqno (1.3)$$

We shall call $L^\infty(G,\Lambda^1(G),\kappa)$ the space of the elements of $\L$ which are bounded $\kappa$-a.e., with the norm of the essential supremum. 

\vskip 1pc

\noindent{\bf The gasket is totally geodesic.} In [16] and [14] it is proven that any two points $x,y\in G$ are joined by a shortest curve with image in $G$; we define $d(x,y)$ as the length of this curve. It is standard that $d$ is a distance on $G$; in [14] it is shown that the two topologies on $G$, the one induced by the immersion in $\R^2$ and the one induced by $d$, coincide. From now on, we shall work on the measured metric space 
$(G,d,\kappa)$.

\vskip 2pc

\centerline{\bf \S 2}
\centerline{\bf Differential forms on the gasket}

\vskip 1pc

\noindent{\bf Differential and codifferential.} First of all, we recall an integration by parts formula from [24]. 

We consider the following three vectors of $\R^2$, which play the r\^ole of the exterior normal to the fractal:  
$$\nu_A=\left(
\matrix{-1
\cr
0}
\right)  ,\qquad
\nu_B=\left(
\matrix{\frac{1}{2}
\cr
\frac{\sqrt{3}}{2}}
\right) ,\qquad
\nu_C=\left(
\matrix{\frac{1}{2}
\cr
-\frac{\sqrt{3}}{2}}
\right)  . $$
If $u\in C^2(\R^2,\R)$ we define 
$$\Delta u={\rm Tr}(E_xD^2 u)  $$
where ${\rm Tr}$ denotes the trace, $D^2u$ is the matrix of second order derivatives and 
$E_x$ is as in (1.1). Let $u\in C^2(\R^2,\R)$ and $v\in C^1(\R^2,\R)$; denoting by 
$\partial_\nu u$ the directional derivative, we recall from [24] (an alternative proof is in [7]) that
$$\ec(\dr u,\dr v)=-\int_G\Delta u\cdot v\dr\kappa+
\sum_{j\in \{ A,B,C \}}v(j)\partial_{\nu_j}u(j)   .  \eqno (2.1)$$
This yields ([8]) a formula for the codifferential of $C^1$ one-forms; namely, if 
$u=a\dr x+b\dr y$ with $a,b\in C^1(\R^2,\R)$ and $E_x$ is as in (1.1), we define 
$$\bar\dr_Gu=-[
(E_x\dr a,\dr x)+(E_x\dr b,\dr y)
]  .  \eqno (2.2)$$
Then, for all $\phi\in C^1(\R^2,\R)$ we have 
$$\ec(u,\dr\phi)=
\int_G(\bar\dr_Gu)\cdot\phi\dr\kappa+
\sum_{j\in \{ A,B,C \}}\phi(j)u(j)(\nu_j)    \eqno (2.3)$$
where $u(j)(v_j)$ means that we apply the one-form $u(j)\in\Lambda^1(\R^2)$ to the vector 
$\nu_j\in\R^2$. 

\vskip 1pc

Now we can define differential and codifferential in the spirit of [9]. 

\noindent{\bf Definitions.} We say that the function $\bar\dr_Gu\in L^2(G,\R,\kappa)$ is the weak codifferential of $u\in L^2(G,\Lambda^1(G),\kappa)$ if, for all 
$\phi\in C^1(\R^2,\R)$, the following integration by parts formula holds.  
$$\ec(u,\dr\phi)=\int_G(\bar\dr_Gu)\cdot\phi\dr\kappa  .   \eqno (2.4)$$
Since $C^1(\R^2,\R)$ is dense in $L^2(G,\R,\kappa)$, the formula above implies easily that the weak codifferential, if it exists, is unique; we call 
$\dc(\bar\dr_G)\subset L^2(G,\Lambda^1(G),\kappa)$ the domain of the codifferential. 

For the vectors $\nu_A$, $\nu_B$, $\nu_C$ defined at the beginning of this section we set 
$$C^1_{abs}(\R^2,\Lambda^1(\R^2))\colon=\{
u\in C^1(\R^2,\Lambda^1(\R^2))\st u(j)(\nu_j)=0\txt{for}
j\in\{ A,B,C \}
\}  .  $$
Now formula (2.3) implies that $C^1_{abs}(\R^2,\Lambda^1(\R^2))\subset\dc(\bar\dr_G)$; moreover, if $u\in C^1_{abs}(\R^2,\Lambda^1(\R^2))$, then its weak codifferential is given by (2.2). The notation for $C^1_{abs}$ is taken from the "absolute" boundary conditions of [22], which generalise Neumann's to $k$-forms. 

We just saw that $C^1_{abs}(\R^2,\Lambda^1(\R^2))\subset\dc(\bar\dr_G)$; since clearly 
$$\{
E_x u(x)\st u\in C^1_{abs}(\R^2,\R)
\}$$
is dense in $L^2(G,\Lambda^1(G),\kappa)$, we get that 
$\dc(\bar\dr_G)$ is dense in $L^2(G,\Lambda^1(\R^2),\kappa)$. 

We say that the one-form $\dr_w v\in L^2(G,\Lambda^1(G),\kappa)$ is the weak differential of 
$v\in L^2(G,\R,\kappa)$ if, for all $u\in\dc(\bar\dr_G)$, the following integration by parts formula holds.  
$$\ec(u,\dr_w v)=\int_G(\bar\dr_Gu)\cdot v\dr\kappa  .  \eqno (2.5)$$
Since we saw above that $\dc(\bar\dr_G)$ is dense in $L^2(G,\Lambda^1(G),\kappa)$, we get as above that the weak differential is unique. 

We call $\dc(\dr_w)\subset L^2(G,\R,\kappa)$ the domain of the weak differential; the definition of the codifferential immediately implies that $C^1(\R^2,\R)\subset\dc(\dr_w)$ and that the weak differential 
$\dr_wu$ satisfies 
$$\dr_w u=E_x\dr u(x)$$
when $u\in C^1(\R^2,\R)$. From now on, we shall often write $\dr u$ instead of $\dr_wu$. 

\vskip 1pc

We recall the proof that the weak differential $\dr$ is closed in $L^2(G,\R,\kappa)$.  

\lem{2.1} Let $\{ v_n \}_{n\ge 1}\subset\dc(\dr)$ be such that 

\noindent 1) $v_n\weak v$ in $L^2(G,\R,\kappa)$ and 

\noindent 2) $\dr v_n\weak V$ in $L^2(G,\Lambda^1(G),\kappa)$. 

Then, $v\in\dc(\dr)$ and $\dr v=V$. 

\proof Since $v_n\in\dc(\dr)$, for all $u\in\dc(\bar\dr_G)$ we have 
$$\ec(\dr v_n,u)=\int_Gv_n(\bar\dr_Gu)\dr\kappa  .  $$
Recall that the inner product in $L^2(G,\Lambda^1(G),\kappa)$ is that of (1.3); now points 1) and 2) above imply respectively the first and second limit below. 
$$\int_G v_n(\bar\dr_Gu)\dr\kappa\tends\int_G u(\bar\dr_Gu)\dr\kappa
\txt{and}
\ec(\dr v_n,u)\tends\ec(V,u)   .  $$
From the last two formulas we get that, for all $u\in\dc(\bar\dr_G)$, 
$$\ec(V,u)=
\int_G v(\bar\dr_Gu)\dr\kappa$$
which implies the thesis by (2.5).

\fin

Now we consider the subspace 
$$\{
(u,\dr u)\st u\in\dc(\dr)
\} \subset 
L^2(G,\R,\kappa)\times L^2(G,\Lambda^1(G),\kappa)   $$
which we endow with the natural inner product   
$$(u,v)_{\dc(\ec^\dr)}\colon=
(u,v)_{L^2(G,\R,\kappa)}+\ec(\dr u,\dr v)  .  \eqno (2.6)$$
The last lemma implies that this subspace is closed for the topology induced by (2.6), which is equivalent to say that the bilinear form 
$$\fun{\ec^\dr}{\dc(\dr)\times\dc(\dr)}{\R}$$
$$\ec^\dr(u,v)=\int_G(\dr u\cdot E_x,\dr v)\dr\kappa$$
is closed in $L^2(G,\kappa)$, i.e. it is a Dirichlet form. Since it is proven in [9] that the operator $\dr$ is local, we get that the Dirichlet form $\ec^\dr$ is local. Defining $\ec^\dr$ as in the formula above obviously implies that $\dc(\ec^\dr)=\dc(\dr)$ and we shall often use this fact in the following. In lemma 2.2 below we are going to show that $C^1(\R^2,\R)$ is dense in $\dc(\dr)$ for the topology of (2.6); in other words, the closure of $\{ (u,\dr u)\st u\in C^1(\R^2,\R) \}$ is the graph of $\ec^\dr$ on 
$\dc(\dr)$. 

\noindent{\bf The Poincar\'e\ inequality.} For starters, we recall a standard notation: if 
$u\in L^1(G,\kappa)$ and $F\subset G$ is a Borel set with $\kappa(F)>0$, we set 
$$u_F=\frac{1}{\kappa(F)}\int_F u\dr\kappa  .  \eqno (2.7)$$
From now on, $B(x,r)$ will denote the open ball of radius $r$ and centre $x\in G$ with respect to the geodesic distance $d$ on $G$. Proposition 4.26 of [14] says that there are $C>0$ and 
$\l\ge 1$ such that for all $u\in\dc(\dr)$, $x\in G$ and $r>0$ we have 
$$\int_{B(x,r)}|u-u_{B(x,r)}|^2\dr\kappa\le
Cr^2\int_{B(x,\l r)}(E_x\dr u,\dr u)\dr\kappa  .  \eqno (2.8)$$

We recall from theorem 2.6 of [14] that on $G$ there is a distance $\rc^\2$ (which is called the resistance distance) such that, if $u\in\dc(\dr)$, the first equality below holds; for the second one, $C$ is the $\rc^\2$-diameter of $G$, which is finite since $G$ is compact and $\rc^\2$ induces the same topology as the geodesic distance $d$.
$$|u(x)-u(y)|^2\le \rc(x,y)\ec(\dr u,\dr u)\le
C^2 \ec(\dr u,\dr u)   .  \eqno (2.9)$$
A consequence of the first inequality above is that the functions $u\in\dc(\ec^\dr)$ are continuous; from the second inequality we get there is $C>0$ such that, for all $u\in\dc(\ec^\dr)$, 
$$||u||_{L^\infty}\le C\left[
\min(|u(A)|,|u(B)|,|u(C)|)+\sqrt{\ec^\dr(u,u)}
\right]  .  $$
This implies that, if $\a,\b,\g\in\R$, then the convex set of the $u\in\dc(\ec^\dr)$ such that $u(A)=\a$, $u(B)=\b$ and $u(C)=\g$ is closed for the topology of (2.6). The formula above also implies that there is $C_1>0$, depending only on $\a$, $\b$ and $\g$, such that, for all such $u$, 
$$\int_G|u|^2\dr\kappa\le C_1[1+\ec(\dr u,\dr u)]   .  \eqno (2.10)$$
A particular case is the subspace $H_0$ of the $\phi\in \dc(\ec^\dr)$ such that 
$\phi(A)=\phi(B)=\phi(C)=0$; we saw above that it is closed and the inequality before (2.10) implies that the norm $||\cdot||_{H_0}$ on $H_0$ defined by  
$$||\phi||_{H_0}^2\colon=
\ec(\dr\phi,\dr\phi)   $$
is equivalent to 
$$||\phi||_{\dc(\ec^\dr)}^2=||\phi||_{L^2}^2+\ec(\dr\phi,\dr\phi)  .  $$
We shall use often the norm $||\cdot||_{H_0}$, which is the one induced by the inner product 
$$(u,v)_{H_0}\colon=\ec(\dr u,\dr v)   .    \eqno (2.11)$$
Since we saw above that $H_0$ is closed in the Hilbert space $\dc(\ec^\dr)$, we get that 
$H_0$ is Hilbert for the inner product of (2.11). 

\noindent{\bf Other properties of the differential.} We prove that test functions are dense.  

\lem{2.2} The space $C^1(\R^2,\R)$ is dense in $\dc(\dr)=\dc(\ec^\dr)$ with the topology of (2.6). 

\proof When we defined the operator $\dr$ we saw that $C^1(\R^2,\R)\subset\dc(\dr)$; now we show that the orthogonal space of $C^1(\R^2,\R)$ in $\dc(\ec^\dr)$ with the natural inner product of (2.6) is reduced to zero. 

Let $u\in\dc(\dr)$ be such that 
$$\int_G u\cdot\phi\dr\kappa+
\int_G(E_x\dr u,\dr\phi)\dr\kappa=0
\txt{for all} \phi\in C^1(\R^2,\R)  .  \eqno (2.12)$$
First of all we assert that, with the notation of (2.7), $u_G=0$. Let us suppose by contradiction that this is not the case; then, we can add a constant to $\phi$ in such a way that the first integral above is zero; since $\dr(\phi+c)=\dr\phi$ we have that 
$$\int_G(E_x\dr u,\dr\phi)\dr\kappa=0
\txt{for all} \phi\in C^1(\R^2,\R)  .  \eqno (2.13)$$
This means that $u$ is harmonic; by [14], this implies that $u$ is affine. Since in the integration by parts formula (2.13) there are no boundary terms, formula (2.1) implies that 
$\partial_{\nu_A}u(A)=\partial_{\nu_B}u(B)=\partial_{\nu_C}u(C)=0$; since $u$ is affine, we get that $u$ is constant; now (2.12) implies that 
$$\int_G u\cdot\phi\dr\kappa=0$$
for all $\phi\in C^1(\R^2,\R)$, i.e. that $u\equiv 0$, which contradicts $u_G\not=0$. 

Thus, if $u$ is orthogonal to $C^1$ we must have that $u_G=0$ and now it suffices to show that the functions $C^1$ with zero average are dense among the functions $u\in\dc(\ec^\dr)$ with $u_G=0$. Formula (2.8) implies that $\ec^\dr(\cdot,\cdot)$ is an equivalent inner product in the subspace of the functions $u$ with $u_G=0$ and we are reduced to show that, if $u$ with $u_G=0$ satisfies (2.13) for all $\phi\in C^1(\R^2,\R)$ with $\phi_G=0$, then $u\equiv 0$. Since again $\dr(\phi+c)=\dr\phi$, we get that $u$ satisfies (2.13) for all $\phi\in C^1(\R^2,\R)$; this implies as above that $u$ is constant; since $u_G=0$, we get that $u\equiv 0$, as we wanted.

\fin

\lem{2.3} 1) Let $u\in\dc(\dr)$ and $v\in\dc(\bar\dr_G)$; let us suppose in addition that 
$u\in L^\infty(G,\R,\kappa)$ and $v\in L^\infty(G,\Lambda^1(G),\kappa)$. Then, 
$u\cdot v\in \dc(\bar\dr_G)$ and 
$$\bar\dr_G(u\cdot v)=
-(E_x\dr u,v)+u\bar\dr_Gv  .  \eqno (2.14)$$
The signs in the formula above are due to the fact that $\bar\dr_G$ is minus the divergence. 

\noindent 2) Let $\psi\in C^1(\R,\R)$ have bounded derivative and let $u\in\dc(\dr)$. Then, $\psi\circ u\in\dc(\dr)$ and 
$$\dr(\psi\circ u)=\psi^\prime\circ u\cdot\dr u  .  \eqno (2.15)$$

\noindent 3) Let $u\in\dc(\dr)$, let $c\in\R$ and let us define 
$$v(x)=\max(c,u(x))  .  $$
Then, $v\in\dc(\dr)=\dc(\ec^\dr)$ and for $\kappa$-a.e. $x\in G$ we have 
$$\dr v(x)=\left\{
\eqalign{
0&\txt{if}u(x)\le c\cr
\dr u(x)&\txt{if}u(x)>c  .
}
\right.   $$

\proof We begin with point 1); note that, by our hypotheses, the right hand side of (2.14) is in $L^2(G,\R,\kappa)$. For starters, we suppose that $u\in C^1(\R^2,\R)$. 

Let $\phi\in C^1(\R^2,\R)$; the first equality below follows from the definition of $\ec$ in (1.1); the second one is obvious and the third one comes from the Leibnitz rule for $C^1$ functions; the last one is the integration by parts formula (2.4), which we can apply because $v\in\dc(\bar\dr_G)$ by assumption. 
$$\ec(uv,\dr\phi)=
\int_G(E_x(uv),\dr\phi)\dr\kappa=
\int_G(E_x v,u\dr\phi)\dr\kappa=$$
$$\int_G(E_x v,\dr(u\phi))\dr\kappa-\int_G\phi(E_xv,\dr u)\dr\kappa=$$
$$\int_G(\bar\dr_G v)u\phi\dr\kappa-\int_G\phi(E_xv,\dr u)\dr\kappa  .  $$
This proves (2.14) when $u\in C^1(\R^2,\R)$; we recall the standard density argument for the general case. By lemma 2.2, $C^1(\R^2,\R)$ is dense in $\dc(\dr)$ with the topology of (2.6); thus, we can find a sequence 
$\{ u_n \}\subset C^1(\R^2,\R)$ such that 
$$u_n\tends u,\txt{and} \dr u_n\tends\dr u
\txt{in}
L^2  .  \eqno (2.16)$$
Since $v,\phi,\dr\phi,E_x\in L^\infty$, we can take limits in the left and right hand sides of 
$$\ec(u_nv,\dr\phi)=
\int_G(\bar\dr_G v)u_n\phi\dr\kappa-\int_G\phi(E_xv,\dr u_n)\dr\kappa    $$
getting (2.14) in the general case. 

For the proof of (2.15) we consider the sequence $\{ u_n \}$ of (2.16); up to taking a subsequence, which we denote with the same index, we can suppose that the convergences of (2.16) are $\kappa$-a.e. and dominated. Since $\psi^\prime$ is bounded, 
this implies that 
$$\psi\circ u_n\tends\psi\circ u
\txt{and}
\dr(\psi\circ u_n)\tends\dr(\psi\circ u)  \eqno (2.17)$$
$\kappa$-a.e. and dominated and thus in $L^2(G,\kappa)$. Since $\psi\circ u_n\in C^1(\R^2,\R)$, the definition of the codifferential implies that for all $\phi\in\dc(\bar\dr_G)$ we have the first equality below, while the second one comes from the chain rule for $C^1$ functions. 
$$\int_G\psi\circ u_n\cdot\bar\dr_G\phi\dr\kappa=
\int_G(E_x\dr(\psi\circ u_n),\phi)\dr\kappa=
\int_G\psi^\prime\circ u_n(E_x\dr u_n,\phi)\dr\kappa  .  $$
Using (2.17) and the fact that $\psi^\prime\circ u_n$ is bounded, we can take limits in the formula above and get  
$$\int_G\psi\circ u\cdot\bar\dr_G\phi\dr\kappa=
\int_G\psi^\prime\circ u(E_x\dr  u,\phi)\dr\kappa    $$
which implies (2.15) by the definition of $\dr$ in (2.5). 

As for point 3), it just says that the operator $\dr$ is local, a fact proven in [9]. Let us recall a different proof, equally standard. We find a sequence $\{\psi_n \}\subset C^1(\R,\R)$ with bounded derivatives such that $\psi_n(s)= c$ when $s\le c$ and $\psi_n(s)\tends s$ when $s>c$; we can also require that $\psi^\prime_n(s)=1$ when $s\ge c+\frac{1}{n}$. We apply point 2) to the composition $\psi_n\circ u$ and then  we take limits.

\fin

\noindent{\bf The codifferential in the relative sense and the Dirichlet problem.} We define the codifferential in the relative sense, i.e we ask that (2.4) holds, but only for test functions which vanish on the boundary of the fractal.

\noindent{\bf Definitions.} Let $v\in L^2(G,\Lambda^1(\R^2),\kappa)$; we say that 
$\bar\dr_G v\in L^2(G,\R,\kappa)$ is the codifferential in the relative sense if formula (2.4) holds for all $\phi\in C^1(\R^2,\R)$ such that $\phi(A)=\phi(B)=\phi(C)=0$. Formula (2.3) implies easily that, if $u\in C^1(\R^2,\Lambda^1(\R^2))$, then it has a codifferential in the relative sense which is given by (2.2). 

Let $f\in L^1(G,\R,\kappa)$, let $g\in L^2(G,\Lambda^1(G),\kappa)$ and let $H_0$ be the space defined before (2.11); we say that $u\in H_0$ is a weak solution of the equation 
$$-\Delta u=f+\bar\dr_G g  \eqno (2.18)$$
if for all $\phi\in C^1(\R^2,\R)$ with $\phi(A)=\phi(B)=\phi(C)=0$ we have that, for the weak differential operator $\dr$ of (2.5), 
$$\ec(\dr u,\dr\phi)=
\int_Gf\cdot\phi\dr\kappa+\int_G(\dr\phi,g)\dr\kappa  .  \eqno (2.19)$$
We have the following. 

\lem{2.4} Let $f\in L^1(G,\R,\kappa)$ and let 
$g\in L^2(G,\Lambda^1(G),\kappa)$; then, equation (2.18) has a unique weak solution 
$u\in H_0$. Moreover, there is $C>0$, independent of $f$ and $g$, such that 
$$||u||_{H_0}\le C\left(
||f||_{L^1}+||g||_{L^2}
\right)   .  $$

\proof Recall that the space $H_0$ is Hilbert for the inner product $(\cdot,\cdot)_{H_0}$ of (2.11). The inequality after (2.9) implies that, if $u\i H_0$, 
$$||u||_{L^\infty}\le C||u||_{H_0}  .  $$
Together with H\"older's inequality this implies that the linear functional 
$$\fun{L}{H_0}{\R}$$
$$\fun{L}{\phi}{
\int_Gf\cdot\phi\dr\kappa+\int_G(\dr\phi,g)\dr\kappa
}  $$
is continuous for $||\cdot||_{H_0}$; actually, by the last two formulas a bound on the operator norm of $L$ is $C(||f||_{L^1}+||g||_{L^2})$. Now the thesis follows from the standard Lax-Milgram argument: the functional $L$ is represented by a unique $u\in H_0$, which satisfies (2.19) for all $\phi\in C^1(\R^2,\R)$ with $\phi(A)=\phi(B)=\phi(C)=0$. 

\fin

\vskip 2pc

\centerline{\bf \S 3}
\centerline{\bf The function spaces}

\vskip 1pc

As in the last section, we consider the harmonic gasket $G$ with the geodesic distance $d$ and Kusuoka's measure $\kappa$. By (2.9), the functions of $\dc(\ec^\dr)$ are $\2$-H\"older with respect to the distance $\rc^\2$; since this distance does not coincide with the geodesic distance $d$, we recall some results from [1] and [5]. 

\noindent{\bf Volume doubling.} We recall from [16] and theorem 4.25 of [14] that $(G,d,\kappa)$ has the volume doubling property: there is $C>0$ such that, for all $x\in G$ and $r>0$, we have 
$$0<\kappa(B(x,2r))\le C\kappa(B(x,r))   .   $$
A consequence is that there are $C,Q>0$ such that, for all $0<r\le R$ we have 
$$\frac{
\kappa(B(x,R))
}{
\kappa(B(x,r))
}  \le \left(
\frac{R}{r}
\right)^Q   .   $$
The number $Q$, which in many formulas takes the r\^ole of the dimension, is computed in [14], where it is shown that 
$$Q=\log_5 15\in (1,2)  .  \eqno (3.1)$$

\noindent{\bf Cheeger's energy.} We denote by ${\rm Lip}(G,\R)$ the space of Lipschitz functions on $G$; if $\fun{u}{G}{\R}$ is Lipschitz, we define the local Lipschitz constant as 
$$Lip(u,x)\colon=\limsup_{r\tends 0}\sup_{y\in B(x,r)}
\frac{|u(x)-u(y)|}{r}   .   $$
We are going to relax this function in $L^2(G,\R,\kappa)$. Let $u\in L^2(G,\R,\kappa)$ and let us suppose that there is $h\in L^2(G,\R,\kappa)$ and a sequence 
$\{ u_n \}\subset {\rm Lip}(G,\R)$ such that the following two points hold. 

\noindent 1) $u_n\tends u$ in $L^2(G,\kappa)$ and 

\noindent 2) $Lip(u_n,\cdot)\weak h$ in $L^2(G,\kappa)$. 

If the set of such functions $h$ is nonempty, then ([4], [13]) it contains a unique element 
$|Du|_\ast$ which is minimal in the following sense: for all $h$ as in points 1) and 2) above, we have 
$$|| |Du|_\ast ||_{L^2(G,\R,\kappa)}\le||h||_{L^2(G,\R,\kappa)}
\txt{and}
|Du|_\ast(x)\le h(x)
\txt{for $\kappa$-a.e} x\in G   .  $$

We define Cheeger's energy as 
$$Ch(u)=\int_G|Du|_\ast^2\dr\kappa  $$
if $|Du|_\ast$ exists, and we set $Ch(u)=+\infty$ otherwise. 

Clearly, the minimality of $|Du|_\ast$ implies that, if $u\in {\rm Lip}(G,\R)$, then 
$$|Du|_\ast(x)\le Lip(u,x)
\txt{for $\kappa$-a.e} x\in G  .  $$

A result of [17]  says that, on the harmonic gasket with the geodesic distance and Kusuoka's measure, $\dc(Ch)=\dc(\ec^\dr)$ and that for all 
$u\in\dc(\ec^\dr)=\dc(\dr)$ we have that 
$$|Du|_\ast^2=
(\dr u(x)E_x,\dr u(x))
\txt{for $\kappa$-a.e} x\in G  .  \eqno (3.2)$$
The last formula and the Poincar\'e\ inequality (2.8) imply that, 
if $u\in\dc(\ec^\dr)=\dc(\dr)$, 
$$\int_{B(x,r)}|u-u_{B(x,r)}|^2\dr\kappa\le
Cr^2\int_{B(x,\l r)}|Du|_\ast^2\dr\kappa(x)  .  \eqno (3.3)$$

\noindent{\bf The Korevaar-Schoen energy.} We recall some definitions and results from [1] and [5]. 

First of all, if $u\in L^2(G,\R,\kappa)$, we define 
$$E_{2,1}(u,r)\colon=
\int_G\dr\kappa(x)\frac{1}{\kappa(B(x,r))}
\int_{B(x,r)}\frac{|f(y)-f(x)|^2}{r^2}\dr\kappa(y)$$
and 
$$||u||_{\bc^{1,2}}^2\colon=
||u||_{L^2(G,\kappa)}^2+\limsup_{r\tends 0}E_{2,1}(u,r)  .   \eqno (3.4)$$
This is known as the Korevaar-Schoen norm of $u$; the space $\kc\sc$ of the functions 
$u$ such that $||u||_{\bc^{1,2}}<+\infty$ is a Banach space for the norm 
$||\cdot||_{\bc^{1,2}}$ of (3.4). 

We saw at the beginning of this section that $(G,d,\kappa)$ has the volume doubling property; since the Poincar\'e\ inequality of (3.3) holds, we are in the hypotheses of theorem 3.1 of [1]; by point (v) of this theorem, $||\cdot||_{\bc^{1,2}}$ is equivalent to the (apparently stronger) norm $||\cdot||_{\bullet}$ defined by 
$$||u||_{\bullet}^2\colon=||u||^2_{L^2(G,\kappa)}+
\sup_{r>0}E_{2,1}(u,r)  .  \eqno (3.5)$$

For the number $Q$ of (3.1) we set 
$$\a\colon=1-\frac{Q}{2}\in\left( 0,\2 \right)     $$
and we recall the standard notation  
$$|u|_{0,\a}\colon=\sup\left\{
\frac{|u(x)-u(y)|}{d(x,y)^\a}\st x\not=y\in G
\right\}     $$
for the H\"older seminorm, and  
$$||u||_{0,\a}\colon=||u||_{\sup}+|u|_{0,\a}    $$
for the H\"older norm. 

Now theorem 3.2 of [5] says that there is $C>0$ such that, for all $u\in\kc\sc$, 
$$||u||_{0.\a}\le C||u||_{\bullet}  .  $$
We saw before (3.5) that $||\cdot||_{\bullet}$ and $||\cdot||_{\bc^{1,2}}$ are equivalent; thus, possibly enlarging $C$, for all $u\in\kc\sc$, 
$$||u||_{0,\a}\le C||u||_{\bc^{1,2}}  .  \eqno (3.6)$$
In [1] it is shown that the $\Gamma$-limit of $E_{2.1}(\cdot,r)$ as $r\tends 0$ is comparable to Cheeger's energy, which implies the second inequality below; for the first inequality we use theorem 3.1 of [1], which says that the $\Gamma$-limit is comparable to $\limsup_{r\tends 0}E_{2,1}(u,r)$. 
$$\limsup_{r\tends 0}E_{2,1}(u,r)\le
C\cdot\Gamma-\lim_{r\tends 0}E_{2,1}(u,r)
\le C_1\cdot Ch(u) .  $$
By (3.2) this implies that $\dc(\dr)\subset{\cal KS}$; together with (3.4), it also implies the second inequality below for all $u\in\dc(\dr)$; the first inequality is (3.6) and the equality at the end comes again from (3.2). 
$$||u||_{0,\a}\le 
C||u||_{\bc^{1,2}}\le
C\left[
||u||_{L^2(G,\kappa)}^2+Ch(u)
\right]^\2   =    C\left[
||u||_{L^2(G,\kappa)}^2+\ec(\dr u,\dr u)
\right]^\2  .  \eqno (3.7)$$
This formula will be essential in the next sections.

\vskip 2pc

\centerline{\bf \S 4}
\centerline{\bf Existence of the minimum } 

\vskip 1pc

\noindent{\bf Invariance under automorphisms.} We recall the classical proof that the functional of (5) is invariant under automorphisms of $H_+$. 

\lem{4.1} Let $\fun{\zeta}{H_+}{H_+}$ be an automorphism and let $\tilde\ec$ be as in (5). Then, for all $u+iv\in\dc(\dr)\oplus i\dc(\dr)$ we have that
$$\tilde\ec(\zeta\circ(u+iv))=\tilde\ec(u+iv)  .  $$

\proof After identifying $\C$ with $\R^2$, we write $\left(
\matrix{
u\cr
v
}
\right)=
\left(
\matrix{
u^1\cr
u^2
}
\right)$, $\zeta=
\left(
\matrix{
\zeta^1\cr
\zeta^2
}
\right)$
and we denote by $g_{i,j}$ the Riemannian tensor of $H_+$. We denote by 
$\left(
\matrix{
x^1\cr
x^2
}
\right)$
the coordinates on $G$ and, for the matrix field $E_x$ of (2) we set 
$E(x^1,x^2)=a^{i,j}(x^1,x^2)$ and adopt the usual convention of summation on repeated indices. With this notation, the definition (5) becomes 
$$\tilde\ec(u^1+iu^2)=
\int_Ga^{i,j}\cdot g_{l,k}\circ(u+iv)\partial_iu^l\partial_ju^k\dr\kappa  .  $$
This implies the last equality below; together with the chain rule, it implies the first one; the second equality follows from the fact that $\zeta$, being an automorphism, preserves the Riemannian tensor of $H_+$. 
$$\tilde\ec(\zeta\circ(u^1+iu^2))=
\int_G a^{i,j}\cdot g_{l,k}\circ\zeta\circ(u^1+iu^2)\cdot
\partial_r\zeta^l|_{u^1+iu^2}\partial_iu^r\partial_s\zeta^k|_{u+iv}\partial_ju^s
\dr\kappa=$$
$$\int_Ga^{i,j}\cdot g_{r,s}\circ(u+iv)\partial_iu^r\partial_ju^s\dr\kappa=
\tilde\ec(u^1+iu^2)  .  $$

\fin

\noindent{\bf The infimum is finite.} We consider three points of $H_+$, not necessarily distinct; however, if they are the vertices of a non-degenerate Euclidean triangle, we ask that they are labeled in the clockwise way, as in figure 1; more precisely, the automorphism bringing $\tilde A$, $\tilde B$ to $i$, $it$ respectively with $t>1$ brings $\tilde C$ to the right hand plane.
$$\tilde A=a_1+ib_1,\qquad 
\tilde B=a_2+ib_2, \qquad
\tilde C=a_3+ib_3  .   \eqno (4.1)$$

We briefly prove that the number $\a$ of (6) is finite: on one side, $\a\ge 0$ since $\tilde\ec\ge 0$. In order to show that $\a<+\infty$, we build 
$\psi_0\in\ac$ such that $\tilde\ec(\psi_0)<+\infty$. 

We consider the canonical map 
$$\fun{J}{\R^2}{\C},\qquad
\fun{J}{(x,y)}{x+iy  .}   $$
Let $\fun{\tilde\psi_0}{\R^2}{\R^2}$ be the unique affine map  such that 
$$\tilde\psi_0(A)=J^{-1}(\tilde A),\qquad 
\tilde\psi_0(B)=J^{-1}(\tilde B),\qquad
\tilde\psi_0(C)=J^{-1}(\tilde C)  .  $$
We define $\psi_0=u_0+iv_0$ as  
$$\psi_0=J\circ\tilde\psi_0  .  $$
This is an affine map from $\R^2$ to $\C$ bringing $A$, $B$, $C$ respectively to $\tilde A$, $\tilde B$, $\tilde C$; in particular, it belongs to the set $\ac$ defined in the introduction. Since $\psi_0$ is affine, $\dr u_0$, $\dr v_0$ are constant one-forms, while $v_0|_G$ is bounded away from $0$ and 
$+\infty$; now (5) implies that $\tilde\ec(\phi_0)<+\infty$. 

\noindent{\bf A bounded minimising sequence.} By the last paragraph, there is a minimising sequence $\{ \phi_n \}$ for (6); since we want $\phi_n(G)$ to lie in a hyperbolically bounded set of $H_+$, independent of $n$, we need the next definitions.   

\vskip 1pc

\noindent{\bf Definitions.} Let $D,E\in H_+$; we denote by 
$$\fun{\eta_{D,E}}{[0,1]}{H_+}$$
the unique hyperbolic geodesic with $\eta_{D,E}(0)=D$ and $\eta_{D,E}(1)=E$. We recall that $\eta_{D,E}$ is an arc of a circle with diameter on the real axis, or a segment on a vertical line. 

We denote by $T(\tilde A,\tilde B,\tilde C)$ the closed hyperbolic triangle triangle with vertices $\tilde A$, $\tilde B$ and $\tilde C$; in other words, $T(\tilde A,\tilde B,\tilde C)$ is the closure of the bounded connected component of 
$$H_+\setminus(
\eta_{\tilde A,\tilde B}([0,1])\cup
\eta_{\tilde B,\tilde C}([0,1])\cup
\eta_{\tilde C,\tilde A}([0,1])
)   .   \eqno (4.2)$$

\vskip 1pc

We have the following. 

\lem{4.2} There is a minimising sequence $\{ \tilde\phi_n \}\subset\ac$ such that 
$$\tilde\phi_n(G)\subset T(\tilde A,\tilde B,\tilde C)
\txt{for all}n\ge 1  .  \eqno (4.3)$$

\proof Since $\a<+\infty$, there is a minimising sequence $\{ \phi_n \}\subset\ac$; we are going to use $\{ \phi_n \}$ to build the minimising sequence $\{ \tilde\phi_n \}$ of the thesis. The idea is the following: up to an automorphism of $H_+$, we can suppose that 
$\tilde A=i$, $\tilde B=ti$ with $t\ge 1$ and $\tilde C$ is in the right hand half plane. We are going to compose $\phi_n$ with a map that is the identity on the right hand half plane and projects the left hand half plane onto the imaginary axis; this will put $\phi_n(G)$ in the right hand half plane, i.e on the right side of $\eta_{\tilde A,\tilde B}([0,1])$. The composition continues to be minimising; indeed, we shall see that its energy is smaller than that of $\phi_n$. Repeating this construction for all sides of the triangle, we get 
$\tilde\phi_n$. 

Now to the rigorous proof. 

First of all, we consider a cyclical permutations of $(\tilde A,\tilde B,\tilde C)$ 
$$(\tilde L,\tilde M,\tilde N)\in\{
(\tilde A,\tilde B,\tilde C) ,
(\tilde B,\tilde C,\tilde A) ,
(\tilde C,\tilde A,\tilde B)  
\}     $$  
and we make the following assertion.

Let $(\tilde L,\tilde M,\tilde N)$ be as in the formula above; let $\{ \phi_n \}$ be minimising and let $\zeta_{\tilde L,\tilde M}$ be the automorphism of $H_+$ such that 
$\zeta_{\tilde L,\tilde M}(\tilde L)=i$ and $\zeta_{\tilde L,\tilde M}(\tilde M)=it$ with $t\ge 1$. Since $\zeta_{\tilde L,\tilde M}$ preserves orientation, our labelling of the vertices implies that 
$\zeta_{\tilde L,\tilde M}(\tilde N)$ is in the closed right hand half plane. We assert that we can build another minimising sequence $\{ \tilde\phi_n \}\subset\ac$ such that 
$\zeta_{\tilde L,\tilde M}\circ\tilde\phi_n(G)$ is contained in the closed right half plane. 

Before proving the assertion, we show that it implies the thesis. Indeed, applying it three times to the three cyclical permutations of $(\tilde A,\tilde B,\tilde C)$, we get a new minimising sequence 
$\{ \tilde\phi_n \}$ such that the three sets 
$$\zeta_{\tilde A,\tilde B}\circ\tilde\phi_n(G),\qquad
\zeta_{\tilde B,\tilde C}\circ\tilde\phi_n(G),\qquad
\zeta_{\tilde C,\tilde A}\circ\tilde\phi_n(G)  $$
are contained in the right hand half plane; by the definitions of the maps 
$\zeta_{\tilde L,\tilde M}$ this is equivalent to (4.3), ending the proof. 

Now we prove the assertion; to fix ideas we suppose that 
$(\tilde L,\tilde M,\tilde N)=(\tilde A,\tilde B,\tilde C)$. By lemma 4.1, for all $n\ge 1$ we have that 
$$\tilde\ec(\zeta_{\tilde A,\tilde B}\circ\phi_n)=\tilde\ec(\phi_n)  .  \eqno (4.4)$$
We consider the function "positive part" 
$$m(x)=\left\{
\eqalign{
0&\txt{if} x\le 0\cr
x&\txt{if} x>0    
}
\right.  $$
and we define
$$\fun{g}{H_+}{H_+}$$
$$g(x+iy)=m(x)+iy   \eqno (4.5)$$
which is the identity on ${\rm Re}(z)>0$ and projects ${\rm Re}(z)\le 0$ orthogonally onto the imaginary axis. We give a name to the real and imaginary parts of 
$g\circ\zeta_{\tilde A,\tilde B}\circ\phi_n$ and $\zeta_{\tilde A,\tilde B}\circ\phi_n$, and we define 
$\psi_n$, $\tilde\phi_n$: 
$$\psi_n\colon= g\circ\zeta_{\tilde A,\tilde B}\circ\phi_n=\tilde a_n+ib_n,\qquad
\zeta_{\tilde A,\tilde B}\circ\phi_n=a_n+ib_n  ,  \qquad
\tilde\phi_n\colon=\zeta_{\tilde A,\tilde B}^{-1}\circ\psi_n  .  \eqno (4.6)$$
Note that the imaginary part of $\psi_n$ coincides with that of 
$\zeta_{\tilde A,\tilde B}\circ\phi_n$ by the definition of $g$ in (4.5). 

The first formula of (4.6) and (4.5) imply that $\psi_n(G)$ is contained in the closed right half plane; by the last formula of (4.6) this implies part of the assertion, namely that 
$\zeta_{\tilde A,\tilde B}\circ\tilde\phi_n(G)$ is contained in the right half plane.

Thus, the proof of the lemma is complete if we show first that $\{ \tilde\phi_n \}\subset\ac$, second that $\{ \tilde\phi_n \}$ is minimising. 

We prove that $\{ \tilde\phi_n \}\subset\ac$. First of all, by lemma 2.3 and (4.6) we have that both components of $\tilde\phi_n$ are in $\dc(\dr)$; it remains to show that  
$$\tilde\phi_n(A)=\tilde A,\qquad 
\tilde\phi_n(B)=\tilde B,\qquad
\tilde\phi_n(C)=\tilde C    .   \eqno (4.7)$$
We saw at the beginning of the proof that 
$$\zeta_{\tilde A,\tilde B}(\tilde A),  \zeta_{\tilde A,\tilde B}(\tilde B), 
\zeta_{\tilde A,\tilde B}(\tilde C)\subset 
\{ z\in H_+\st{\rm Re}(z)\ge 0 \}  .  $$
Since by (4.5) $g$ is the identity on the right half plane, the first equality of (4.6) implies that $\psi_n$ brings $A$, $B$, $C$ to $\zeta_{\tilde A,\tilde B}(\tilde A)$, $\zeta_{\tilde A,\tilde B}(\tilde B)$, 
$\zeta_{\tilde A,\tilde B}(\tilde C)$ respectively; by the last one of (4.6), $\tilde\phi_n$ brings $A$, $B$, 
$C$ to $\tilde A$, $\tilde B$, $\tilde C$ respectively, proving (4.7). 

Lastly, we show that $\{ \tilde\phi_n \}$ is minimising. The first and last equalities below come from (4.6) and the definition of $\tilde\ec$ in (5), the second equality follows from  point 3) of lemma 2.3; recall that, by (4.5), $\tilde a_n=m\circ a_n$. The inequality follows since $m^\prime(x)$ is either zero or one. 
$$\tilde\ec(\psi_n)=
\int_G\frac{1}{b_n^2}[
(E_x\dr\tilde a_n,\dr\tilde a_n)+(E_x\dr b_n,\dr b_n)
]\dr\kappa=$$
$$\int_G\frac{1}{b_n^2}[(m^\prime\circ a_n(x))^2\cdot
(E_x\dr a_n,\dr a_n)+(E_x\dr b_n,\dr b_n)
]\dr\kappa\le$$
$$\int_G\frac{1}{b_n^2}[
(E_x\dr a_n,\dr a_n)+(E_x\dr b_n,\dr b_n)
]\dr\kappa  =\tilde\ec(\zeta_{\tilde A,\tilde B}\circ\phi_n)  .  $$
The first equality below is the definition of $\tilde\phi_n$ in (4.6), the second one follows from lemma 4.1 since $\zeta_{\tilde A,\tilde B}$ is an automorphism. 
$$\tilde\ec(\tilde\phi_n)=
\tilde\ec(\zeta_{\tilde A,\tilde B}^{-1}\circ\psi_n)=\tilde\ec(\psi_n)  .  $$
The last  two formulas together with (4.4) imply that 
$$\tilde\ec(\tilde\phi_n)\le\tilde\ec(\phi_n)  .   $$
Since $\{ \phi_n \}$ is minimising,  $\{ \tilde\phi_n \}$ is minimising too.

\fin

\noindent{\bf A bound on the minimal.} The technique of the last lemma yields the following bound, which will be useful in the next section. 

\lem{4.3} Let $\phi\in\ac$ be minimal; then, with the notation of (4.2), 
$$\phi(G)\subset T(\tilde A,\tilde B,\tilde C)   .  \eqno (4.8)$$

\proof As we saw in the last lemma it suffices to show that

\noindent $i$) if $\zeta$ is the automorphism of $H_+$ bringing $\tilde A$ to $i$ and 
$\tilde B$ to $ti$ with $t\ge 1$, or 

\noindent $ii$) if $\zeta$ is the automorphism of $H_+$ bringing $\tilde B$ to $i$ and 
$\tilde C$ to $ti$ with $t\ge 1$, or 

\noindent $iii$) if $\zeta$ is the automorphism of $H_+$ bringing $\tilde C$ to $i$ and 
$\tilde A$ to $ti$ with $t\ge 1$, 

then $\zeta\circ\phi(G)$ is contained in the right half plane. 

We prove case $i$), since the other ones are analogous. 

We begin to note that, since the functions in $\ac$ are continuous, 
$$O\colon=\{
x\in G\st {\rm Re}\zeta\circ\phi(x)<0
\}   \eqno (4.9)$$
is an open set of $G$. We suppose by contradiction that this set is non-empty. 

We set $\zeta\circ\phi=u+iv$ and define $\tilde\phi$ as in (4.6) by $\zeta\circ\tilde\phi=g\circ(u+iv)$ where $g$ is as in (4.5); formula (5) implies that 
$$\tilde\ec(\zeta\circ\phi)=\int_{G\setminus O}
\frac{1}{v^2}[
(\dr uE_x,\dr u)+(\dr vE_x,\dr v)
]\dr\kappa+$$
$$\int_{ O}
\frac{1}{v^2}[
(\dr uE_x,\dr u)+(\dr vE_x,\dr v)
]\dr\kappa   .    $$
By (4.5), $\zeta\circ\tilde\phi=\zeta\circ\phi$ on $G\setminus O$; since the first component of 
$\zeta\circ\tilde\phi$ is zero on $\oc$, we get the equality below. 
$$\tilde\ec(\zeta\circ\tilde\phi)=\int_{G\setminus O}
\frac{1}{v^2}[(\dr uE_x,\dr u)+(\dr vE_x,\dr v)]\dr\kappa+
\int_{O}
\frac{1}{v^2}(\dr vE_x,\dr v)\dr\kappa  .    $$
We assert that 
$$\int_O\frac{1}{v^2}(E_x\dr u,\dr u)\dr\kappa>0  .  $$
Indeed, if the integral in the last formula were zero, then the Poincar\'e\ inequality (2.8) would imply that $u$ is constant on every connected component of $O$. If we show that this constant is zero, we have a contradiction with (4.9). Let $\uc$ be a connected component of $\oc$ and let 
$u(\uc)=\{ x_0 \}$ with ${\rm Re}(x_0)<0$; since $\uc$ is the pre-image of a point it is also closed, and since $G$ is connected we get that $\uc=G$, contradicting the fact that 
${\rm Re}(\tilde C)>0$.

We see as in the last lemma that $\tilde\phi\in\ac$; lemma 4.1 yields the first and last equalities below, the strict inequality comes from the three formulas above. 
$$\tilde\ec(\tilde\phi)=\tilde\ec(\zeta\circ\tilde\phi)<
\tilde\ec(\zeta\circ\phi)=\tilde\ec(\phi)  .  $$
This contradicts the minimality of $\phi$ and we are done. 

\fin 

\noindent{\bf Convergence.} Now we can prove point 1) of theorem 1. 

\prop{4.4} There is a minimising sequence $\{ \phi_n \}\subset\ac$ which, up to subsequences, converges uniformly to a map $\phi\in\ac$ with $\tilde\ec(\phi)=\a$. By the definition in the introduction, the map $\phi$ is harmonic. 

\proof We saw above that there is a minimising sequence $\phi_n=u_n+iv_n$; by lemma 4.2, we can suppose that $\{ \phi_n \}$ satisfies (4.3). It is immediate that there is $C\ge 1$ such that $z\in T(\tilde A,\tilde B,\tilde C)$ implies 
$$\frac{1}{C}\le{\rm Im} z\le C  .  $$
Together with (4.3) this implies that, for all $n\ge 1$, 
$$\frac{1}{C}\le v_n(x)\le C \txt{for $\kappa$-a.e.} x\in G  .  \eqno (4.10)$$
By the definition of $\tilde\ec$ in (5) this implies that 
$$\ec(\dr u_n,\dr u_n)+\ec(\dr v_n,\dr v_n)\le C^2\tilde\ec(\phi_n)  .  $$
Since $\{ \phi_n \}$ is minimising, $\tilde\ec(\phi_n)$ tends to the finite number $\a$; using the formula above we get that there is $C_1>0$ such that, for all 
$n\ge 1$, 
$$\ec(\dr u_n,\dr u_n)+\ec(\dr v_n,\dr v_n)\le C_1   .    \eqno (4.11)$$
Since $u_n\in\ac$, we can apply the Poincar\'e\ inequality (2.10) and get that that there is $C_2>0$ such that, for all $n\ge 1$, 
$$||u_n||_{L^2(G,\R)}+||v_n||_{L^2(G,\R)}\le C_2  .  $$
By (3.7) the last two formulas imply that there is $C_3>0$ such that for all $n\ge 1$ we have 
$$||u_n||_{0,\a}+||v_n||_{0,\a}\le C_3ß  .  \eqno (4.12)$$
By (4.12) and the Arzel\`a-Ascoli theorem we can extract a subsequence (which with we denote with the same index)  
$$\{ \phi_n=u_n+iv_n \}\subset\ac  $$
such that, endowing $C(G,H_+)$ with the $\sup$ norm, 
$$u_n\tends u\txt{and}v_n\tends v\txt{in}C(G,H_+)  .    \eqno (4.13)$$
By (4.11) and weak compactness we can further refine to have 
$$\dr u_n\weak U\txt{and}\dr v_n\weak V\txt{in}L^2(G,\Lambda^1(G),\kappa)  .  $$
By lemma 2.1, the last two formulas imply that that $u,v\in\dc(\dr)$ and that 
$U=\dr u$, $V=\dr v$; thus, the last formula becomes 
$$\dr u_n\weak \dr u\txt{and}\dr v_n\weak \dr v\txt{in}L^2(G,\Lambda^1(\R^2),\kappa)  .  
\eqno (4.14)$$
We just saw that that $u,v\in\dc(\ec^\dr)=\dc(\dr)$; the function $\phi=u+iv$ is our candidate minimum. Since $\phi_n\in\ac$, (4.13) implies that $\phi(A)=\tilde A$, $\phi(B)=\tilde B$, 
$\phi(C)=\tilde C$, yielding that $\phi\in\ac$. The proof is complete if we show that $\phi$ is minimal.

We note that, by (4.10), (4.11) and (4.13)  
$$\lim_{n\tends+\infty}\int_G
\left\vert
\frac{1}{v_n^2}-\frac{1}{v^2}
\right\vert  [
(\dr u_nE_x,\dr u_n)+(\dr v_nE_x,\dr v_n)
]\dr\kappa =0  .  \eqno (4.15)$$
Let now $v\in\ac$ be fixed and let it satisfy $v\ge \frac{1}{C}$ as in (4.10); it is immediate that the convex functional 
$$\fun{}{L^2(G,\Lambda^1(G),\kappa)\times L^2(G,\Lambda^1(G),\kappa)}{\R}$$
$$\fun{}{(\psi,\eta)}{
\int_G\frac{1}{v^2}[
( \psi E_x,\psi )+(\eta E_x,\eta )
]\dr\kappa
}  $$
is continuous for the strong topology of 
$L^2(G,\Lambda^1(G),\kappa)\times L^2(G,\Lambda^1(G),\kappa)$, and thus it is lower semicontinuous for the weak one. 

The first equality below follows by the definition of $\tilde\ec$, the second one by adding and subtracting; the inequality follows from (4.14), (4.15) and weak lower-semicontinuity; the equality at the end is again the definition of $\tilde\ec$. 
$$\lim_{n\tends+\infty}\tilde\ec(\phi_n)=
\lim_{n\tends+\infty}\int_G\frac{1}{v_n^2}\left[
(\dr u_nE_x,\dr u_n)+(\dr v_nE_x,\dr v_n)
\right]\dr\kappa=$$
$$\lim_{n\tends+\infty}\int_G\frac{1}{v^2}\left[
(\dr u_nE_x,\dr u_n)+(\dr v_nE_x,\dr v_n)
\right]\dr\kappa+$$
$$\lim_{n\tends+\infty}\int_G\left(
\frac{1}{v_n^2}-\frac{1}{v^2}
\right)
\left[
(\dr u_nE_x,\dr u_n)+(\dr v_nE_x,\dr v_n)
\right]\dr\kappa  \ge$$
$$\int_G\frac{1}{v^2}\left[
(\dr uE_x,\dr u)+(\dr vE_x,\dr v)
\right]\dr\kappa=\tilde\ec(\phi)  .  $$
We saw above that $\phi\in\ac$; sine $\{ \phi_n \}$ is minimising, the last formula implies that 
$\tilde\ec(\phi)\le\a$, ending the proof. 

\fin

\vskip 2pc

\centerline{\bf \S 5}
\centerline{\bf Uniqueness} 

\vskip 1pc

In this section we prove point 2) of theorem 1; we start with the following lemma.

\lem{5.1} There are $r_0>0$ and $C>0$ such that the following holds. Let $r\in(0,r_0)$ and let 
$$\tilde A,\tilde B,\tilde C\in B(i,r)  ,  \eqno (5.1)$$
where the ball is the Euclidean one. Let $\phi\in\ac$ be minimal; then, 
$$\phi(G)\subset B(i,r)  .  \eqno (5.2)$$

\proof Since $\phi$ is minimal in $\ac$, formula (4.8) holds. Recall that the ball $B(i,r)$ is also a hyperbolic ball, albeit with a different center and radius; since hyperbolic balls are geodesically convex, the geodesics connecting $\tilde A$ with $\tilde B$, $\tilde B$ with $\tilde C$ and $\tilde C$ with $\tilde A$ are contained in $B(i,r)$. As a consequence,  the hyperbolic triangle 
$T(\tilde A,\tilde B,\tilde C)$ satisfies 
$$T(\tilde A,\tilde B,\tilde C)\subset B(i,r)  .  $$
Since $\phi(G)\subset T(\tilde A,\tilde B,\tilde C)$ by (4.8), we get (5.2). 

\fin

Let us call $\tilde r(r)$ the hyperbolic radius of $B(i,r)$; if $B$ is a ball of hyperbolic radius 
$\tilde r(r)$, there is an automorphism of $H_+$ bringing it into $B(i,r)$. Thus, the next proposition implies point 2) of theorem 1. 

\prop{5.2} There is $r_0>0$ such that, if $r\in(0,r_0)$ and (5.1) holds, then $\tilde\ec$ has a unique minimiser $\phi_{\tilde A,\tilde B,\tilde C}$. The function 
$$\fun{\Phi}{B(i,r)^3}{\dc(\ec^\dr)\oplus i\dc(\ec^\dr)} ,\qquad
\fun{\Phi}{(\tilde A,\tilde B,\tilde C)}{\phi_{\tilde A,\tilde B,\tilde C}}  $$
is of class $C^1$. 

\proof By lemma 5.1 and proposition 4.4 it suffices to prove that there is a unique critical point $\phi$ of $\tilde\ec$ which satisfies (5.2) and depends differentiably on $\tilde A$, $\tilde B$, 
$\tilde C$. In order to find such a critical point, in step 1 below we perform a blow-up and recast the Euler-Lagrange equation of $\tilde\ec$ in the form (5.8) below; in step 2, we see that (5.8) means that 
$\phi_{\tilde A,\tilde B,\tilde C}$ is the zero of the nonlinear operator of (5.12), to which we apply the implicit function theorem in step 3. 

\noindent {\bf Step 1: blow up.} We write the function $\phi\in\ac$ in the form 
$\phi= u+i (v+1)$.  By (5), the energy of $\phi$ is
$$\tilde\ec(\phi)=
\int_G\frac{1}{(1+v)^2}[
(\dr u E_x,\dr u)+(\dr v E_x,\dr v)
]\dr\kappa  .  $$
We set 
$$a_r(u,v)=\frac{-2r}{(1+rv)^3}[
(\dr u E_x,\dr u)+(\dr v E_x,\dr v)
]  .  $$
Classical results on the differentiation of functionals ([3]) show that the critical points of $\tilde\ec$ satisfy 
$$\int_G \left\{
\frac{1}{(1+v)^2}[
(\dr u E_x,\dr \phi_1)+(\dr v E_x,\dr \phi_2)]+
a_1(u,v)\phi_2
\right\}\dr\kappa  =0  $$
for all $\phi_1,\phi_2\in\dc(\dr)$ with $\phi_i(A)=\phi_i(B)=\phi_i(C)=0$. By the definition of the codifferential in the relative sense this implies that 
$$\left\{
\eqalign{
\bar\dr_G\frac{\dr u}{(1+v)^2}&=0\cr
\bar\dr_G\frac{\dr v}{(1+v)^2}+a_1(u,v)&=0\cr
\phi=[u+i(v+1)]&\in\ac  .    
}
\right.  $$

In order to make a blow-up around $i$ we set 
$$\phi_r=\frac{\phi-i}{r}+i=u_r+i(1+v_r),\txt{or equivalently}
\phi=u+i(v+1)=ru_r+(1+rv_r)i  .  \eqno (5.3)$$

Formula (5.3) implies the first equality below, while the second one comes from the definition of $\tilde\ec$ in (5). 
$$\tilde\ec(\phi)=
\tilde\ec(ru_r+i(1+rv_r))=$$
$$r^2\int_G\frac{1}{(1+rv_r)^2}[
(\dr u_r E_x,\dr u_r)+(\dr v_r E_x,\dr v_r)
]\dr\kappa  .   $$
Using again the definition of the codifferential in the relative sense, the Euler-Lagrange equation of this functional is given by the first two equations below; for the last three equations, which are the boundary conditions, we use (5.3) and the fact that $\phi(A)=\tilde A$, 
$\phi(B)=\tilde B$, $\phi(C)=\tilde C$.  
$$\left\{
\eqalign{
\bar\dr_G\left(
\frac{\dr u_r}{(1+rv_r)^2}
\right)  &=0\cr
\bar\dr_G\left(
\frac{\dr u_r}{(1+rv_r)^2}
\right)+a_r(u_r,v_r) &=0\cr
(u_r+iv_r)(A)=\frac{1}{r}(\tilde A-i)\cr
(u_r+iv_r)(B)=\frac{1}{r}(\tilde B-i)\cr
(u_r+iv_r)(C)=\frac{1}{r}(\tilde C-i)  .
}
\right.$$
Adding and subtracting, we re-write it in the following way; the last condition is implied by (5.2) and the definition of $u_r$, $v_r$ in (5.3); the ball is the Euclidean one with center the origin and radius 1. 
$$\left\{
\eqalign{
-\Delta u_r+\bar\dr_G\left[
\dr u_r\left(\frac{1}{(1+rv_r)^2}-1\right)
\right]  &=0\cr
-\Delta v_r+\bar\dr_G\left[
\dr v_r\left( \frac{1}{(1+rv_r)^2}-1 \right)
\right] +a_r(u_r,v_r)&=0\cr
(u_r+iv_r)(A)=\frac{1}{r}(\tilde A-i)\cr
(u_r+iv_r)(B)=\frac{1}{r}(\tilde B-i)\cr
(u_r+iv_r)(A)=\frac{1}{r}(\tilde C-i)  \cr
(u_r+iv_r)(G)\subset B(0,1)  .  
}
\right.  \eqno (5.4)$$
Now we consider the affine map $\psi_0$ we defined after (4.1); we make a blow up as in (5.3), setting 
$$\frac{1}{r}[\psi_0(z)-i]=\tilde u_r(z)+i\tilde v_r(z)  .  \eqno (5.5)
$$
Since $\psi_0$ is affine, we have that $\tilde u_r$ and $\tilde v_r$ are affine too and thus by (5.1) 
$\tilde u_r+i\tilde v_r$ has image inside the Euclidean ball $B(0,1)$, which implies the first inequality below; since $\psi_0$ is affine, the first inequality implies the second one for some $C_1>0$ independent of 
$\tilde A,\tilde B,\tilde C\in B(i,r)$. 
$$\sup_{z\in G}|\tilde u_r(z)+i\tilde v_r(z)|\le 1
\txt{and}
\sup_{z\in G}|\dr \tilde u_r(z)+i\dr\tilde v_r(z)|\le C_1    .  \eqno (5.6)$$
In the first equality below, we write $u_r+iv_r$ as a perturbation of the affine map 
$\frac{1}{r}(\psi_0-i)$; the second equality comes from (5.5). 
$$(u_r+iv_r)(z)=\frac{1}{r}[\psi_0(z)-i]+U(z)+iV(z)=\tilde u_r(z)+i\tilde v_r(z)
+U(z)+iV(z)  .    $$
In this way we have that  
$$u_r=\tilde u_r+U,\qquad v_r=\tilde v_r+V  ,  $$
where $\tilde u_r$, $\tilde v_r$ are affine and depend on the points $\tilde A$, $\tilde B$, 
$\tilde C$. 

Since both $\phi$ and $\psi_0$ bring $A$, $B$, $C$ into $\tilde A$, $\tilde B$, $\tilde C$ respectively, (5.3) and (5.5) imply that $(u_r+iv_r)(j)=(\tilde u_r+i\tilde v_r)(j)$ for 
$j\in\{ A,B,C \}$; by the last formula this implies that $U+iV$ vanishes on $A$, $B$ and 
$C$, or that $U,V\in H_0$, where $H_0$ is the space we defined before (2.11). At the beginning of section 2 we saw that the differential of a constant and the Laplacian of an affine function are both zero; thus, the first two equations of (5.4) become the first two equations below, while the third one follows from (5.6) and the last one of (5.4). 
$$\left\{
\eqalign{
-\Delta U+\bar\dr_G\left[
\dr (\tilde u_r+U)\left(\frac{1}{(1+r(\tilde v_r+V))^2}-1\right)
\right]  &=0\cr
-\Delta V+\bar\dr_G\left[
\dr (\tilde v_r+V)\left( \frac{1}{(1+r(\tilde v_r+V))^2}-1 \right)
\right]+a_r(\tilde u_r+U,\tilde v_r+V) &=0\cr
(U+iV)(G)\subset B(0,2),\qquad U,V\in H_0  .
}
\right.  \eqno (5.7)$$
By lemma 2.4 we can apply $(-\Delta)^{-1}$ to the first two equations of (5.7), getting 
$$\left\{
\eqalign{
U+(-\Delta)^{-1}\bar\dr_G\left[
\dr (\tilde u_r+U)\left(\frac{1}{(1+r(\tilde v_r+V))^2}-1\right)
\right]  &=0\cr
V+(-\Delta^{-1})\left\{\bar\dr_G\left[
\dr (\tilde v_r+V)\left( \frac{1}{(1+r(\tilde v_r+V))^2}-1 \right)
\right]
+a_r(\tilde u_r+U,\tilde v_r+V)
\right\} &=0\cr
(U+iV)(G)\subset B(0,2),\qquad U,V\in H_0  .
}
\right.  \eqno (5.8)$$
Note that $U+iV\equiv 0$ is a solution of (5.8) when $r=0$; in the next two steps we shall see that this solution survives when $r>0$ is small. 

\noindent{\bf Step 2. The nonlinear operator.} The solutions of (5.8) are the zeroes of a nonlinear operator; in this step we write it down explicitly and verify the hypotheses of the implicit function theorem. 

We group the nonlinear terms of (5.8) in the function $F$ defined below; it depends on the parameters $\tilde A$, $\tilde B$, $\tilde C$ through the affine function $\tilde u_r+i\tilde v_r$; as for the denominators, we can choose $\e$ so small that $\e|\tilde v_r+V|\le \2$ $\kappa$-a.e. for all $U+iV$ satisfying the last one of (5.8) and we take $r\in(-\e,\e)$. 
$$\fun{F}{(-\e,\e)\times B(0,1)^3\times(H_0\oplus i H_0)}{H_0\oplus i H_0}$$
$$F(r,\tilde A,\tilde B,\tilde C,U+iV)=
(-\Delta)^{-1}\bar\dr_G\left[
\dr (\tilde u_r+U)\left(\frac{1}{(1+r(\tilde v_r+V))^2}-1\right)
\right]+$$
$$i (-\Delta^{-1})\left\{\bar\dr_G\left[
\dr (\tilde v_r+V)\left( \frac{1}{(1+r(\tilde v_r+V))^2}-1 \right)
\right]
+a_r(\tilde u_r+U,\tilde v_r+V)
\right\}  .  \eqno (5.9)$$
We  prove that this operator lands in $H_0\oplus i H_0$. Since $U\in\dc(\ec^\dr)=\dc(\dr)$, we have that $\dr U\in L^2$, while $\dr\tilde u_r\in L^2$ since $\tilde u$ is affine and we saw in section 2 that the differential of a $C^1$ function is the projection of the standard differential on the cotangent space. 

On the other side, (5.6) and the last one of (5.8) imply that, for all $r\in(-\e,\e)$, 
$$\left\vert\left\vert
\left(\frac{1}{(1+r(\tilde v_r+V))^2}-1\right)
\right\vert\right\vert_{L^\infty}\le C_3r   .    $$
By the second inequality of (5.6) this implies that 
$$\left\vert\left\vert
\left(\frac{1}{(1+r(\tilde v_r+V))^2}-1\right)\dr (\tilde u_r+U)
\right\vert\right\vert_{L^2}\le
C_4r(1+||\dr U||_{L^2}),$$
$$\left\vert\left\vert
\left(\frac{1}{(1+r(\tilde v_r+V))^2}-1\right)\dr (\tilde v_r+V)
\right\vert\right\vert_{L^2}\le
C_4r(1+||\dr V||_{L^2})  .  $$
By the definition of $a_r$ at the beginning of step 1 we get simlarly 
$$||a_r(\tilde u_r+U,\tilde v_r+V)||_{L^1}\le
C_4r\left(
1+||\dr U||^2_{L^2}+||\dr V||^2_{L^2}
\right)^\2.    $$
Now lemma 2.4 implies that, for some absolute constant $C_5>0$, 
$$\left\{
\eqalign{
\left\vert\left\vert
(-\Delta)^{-1}\bar\dr_G\left[
\dr (\tilde u_r+U)\left(\frac{1}{(\tilde v_r+rV)^2}-1\right)
\right]
\right\vert\right\vert_{H_0}&\le C_5 r(1+||\dr U||_{L^2})\cr
\left\vert\left\vert
(-\Delta)^{-1}\bar\dr_G\left[
\dr (\tilde v_r+V)\left(\frac{1}{(\tilde v_r+rV)^2}-1\right)
\right]
\right\vert\right\vert_{H_0}&\le C_5 r(1+||\dr V||_{L^2})\cr
||(-\Delta)^{-1}a_r(\tilde u_r+U,\tilde v_r+V)||_{H_0}  &\le
C_5r(1+||\dr U||^2_{L^2}+||\dr V||^2_{L^2})^\2  .
}
\right.    \eqno (5.10)$$
This implies that $F$ lands in a ball of $H_0\oplus i H_0$ which is small with $r$; actually, if $U$, $V$ are in a fixed ball of $H_0$, then (forgetting to write some arguments of $F$) $F(U+iV)$ is in a ball of $H_0\oplus iH_0$ with radius smaller that $C_6r$. 

We assert that $F$ is differentiable in the variables $(U,V)$ and that the differential depends continuously on all variables. Let us consider for instance the real part of $F$; forgetting to write some arguments of $F$ on the left hand side, we have 
$${\rm Re}F(U+iV)=
(-\Delta)^{-1}\bar\dr_G\left[
\dr (\tilde u_r+U)\left(\frac{1}{(1+r(\tilde v_r+V))^2}-1\right)
\right]  .  $$
This operator is affine in $U$ and thus we have that 
$$\partial_{U}{\rm Re}F(U+iV)(h)=
(-\Delta)^{-1}\bar\dr_G\left[
\dr h\left(\frac{1}{(1+r(\tilde v_r+V))^2}-1\right)
\right]  .  \eqno (5.11)$$
This function does not depend on $U$; we show that it is continuous in $r$, $\tilde A$, 
$\tilde B$, $\tilde C$ and $V$. Indeed, the map 
$$\fun{R}{(r,\tilde A,\tilde B,\tilde C,V)}{
\frac{1}{(1+r(\tilde v_r+V))^2}-1
}$$
is continuous from $(-\e,\e)\times B(i,r)^3\times H_0$ to $L^\infty$ by (3.7) and the definition of 
$\tilde v_r$. Now we consider the map 
$$\fun{S}{L^\infty}{\L(H_0,L^2)},\qquad S(W)(h)=W\dr h.   $$
It is easy to see that $S$ is linear and continuous; by lemma 2.4 this implies that the operator 
$$\fun{T}{L^\infty}{\L(H_0,H_0)}$$
$$T(W)(h)=(-\Delta)^{-1}\bar\dr_G(W\dr h)=(-\Delta)^{-1}\bar\dr_G(S(W)h)     $$
is linear and continuous. Since by (5.11)
$$\partial_{U}{\rm Re}F(U+iV)(h)=T\circ R(V)(h)  ,  $$
we get that the partial derivative above is continuous from 
$(-\e,\e)\times B(i,r)^3\times H_0$ to $\L(H_0,H_0)$. 

For the derivative with respect to $V$, standard results on Nemitsky operators (see for instance [3]) show that 
$$\partial_{V}{\rm Re}F(U+iV)(h)=
(-\Delta)^{-1}\bar\dr_G\left[
\dr (\tilde u_r+U)\frac{-2rh}{(1+r(\tilde v_r+V))^3} 
\right]  .  $$
We see as above that this map is continuous from 
$(-\e,\e)\times B(i,r)^3\times H_0\oplus iH_0$ to $\L(H_0,H_0)$. 

The same argument implies that ${\rm Re}F$ is $C^1$ in the variables $r,\tilde A,\tilde B,\tilde C$; we forego the analogous argument for the imaginary part.

Moreover, (5.9) implies that $D_{U,V}F(0,\tilde A,\tilde B,\tilde C,U+iV)$ is the zero operator. Now we define 
$$\fun{G}{(-\e,\e)\times B(0,1)^3\times (H_0\oplus iH_0)}{H_0\oplus iH_0}  $$
$$G(r,\tilde A,\tilde B,\tilde C, U+iV)\colon=
(U+iV)+F(r,\tilde A,\tilde B,\tilde C,U+iV)  .  $$

The map $G$ is $C^1$ in $(-\e,\e)\times B(0,1)^3\times(H_0\oplus iH_0)$; moreover, since 
$D_{U,V}F(0,\tilde A,\tilde B,\tilde C,U+iV)$ is the zero operator, the matrix 
$$\left(
\matrix{
\partial_U G(0,\cdot)&0\cr
0&\partial_VG(0,\cdot)
}
\right)  $$
is the identity from $H_0\oplus iH_0$ to itself and thus it is invertible. 

Equation (5.8) is equivalent to 
$$G(r,\tilde A,\tilde B,\tilde C, U+iV)=0    \eqno (5.12)$$
which has the unique solution $U\equiv V\equiv 0$ when $r=0$.

\noindent{\bf Step 3. The implicit function theorem.} By the last step we can apply the implicit function theorem to (5.12); we get that there is $\e>0$ with the following property. If 
$r\in[0,\e)$ and $\tilde A,\tilde B,\tilde C\in B(0,1)$, then equation (5.12) has a unique solution 
$U_{r,\tilde A,\tilde B,\tilde C}+iV_{r,\tilde A,\tilde B,\tilde C}$; such a solution depends differentiably on $r,\tilde A,\tilde B,\tilde C$ and is zero when $r=0$; in particular, the last one of (5.8) holds. In other words, we are finding a branch of solutions which coincide with the affine $\tilde u_r+i\tilde v_r$ when $r=0$; the dependence on all parameters is $C^1$. 



\fin

\vskip 2pc

\centerline{\bf \S 6}
\centerline{\bf Geodesics}
\vskip 1pc

One of the properties of harmonic immersions is that they are totally geodesic, i.e. they bring geodesics into geodesics. Of course, this is a local property: the image of a minimising curve satisfies the geodesic equation (i.e. it is minimising on small time intervals) but it may not be the shortest path between its endpoints. 

The problem is that at this point we only know that the harmonic map $\phi$ is H\"older, which does not even imply that the image of a geodesic of $G$ has finite length. In this section we are going to show that the image of a geodesic on $G$ has finite length and is a geodesic on $\phi(G)$. 

In order to show that the harmonic immersions $\phi$ is totally geodesic, we recall from [14] and [16] some facts about geodesics on the harmonic Sierpinski gasket $G$ of section 1. 

First of all, there are six geodesics 
$$\fun{\g_{r,s}}{[0,1]}{G},\quad r\not=s\in\{ A,B,C \}  \eqno (6.1)$$
such that 
$$\g_{r,s}(0)=r,\qquad \g_{r,s}(1)=s\txt{for}r\not=s\in\{ A,B,C \}  .  $$
Obviously, if we parametrise these geodesics proportional to arc-length, $\g_{r,s}(t)=\g_{s,r}(1-t)$. 

These geodesics are $C^1$ but not $C^2$ and they are the "outer boundary" of the gasket in the following sense: if $U$ is the unbounded connected component of $\R^2\setminus G$, then 
$$\partial U=\g_{A,B}([0,1])\cup\g_{B,C}([0,1])\cup\g_{A,C}([0,1])   .   $$
On smaller cells we have that the curve connecting, say, $F_{j_0\dots j_{l-1}}(A)$ with 
$F_{j_0\dots j_{l-1}}(C)$, which has minimal length and has image in the cell 
$F_{j_0\dots j_{l-1}}(G)$ is $F_{j_0\dots j_{l-1}}\circ\g_{A,C}$. 

Now the harmonic gasket is totally geodesic: for any two points in $G$ there is a geodesic connecting them. Moreover, this geodesic is a countable or finite chain of geodesics of the type $F_{j_0\dots j_{l-1}}\circ\g_{r,s}$ with $r\not=s\in\{ A,B,C \}$. The chain is finite if and only if the geodesic connects two vertices of the pre-fractal, say 
$\g(0)=F_{j_0\dots j_{l-1}}(A)$ and $\g(1)=F_{i_0\dots i_{l-1}}(B)$.

This prompts us to state the following proposition; its proof will occupy the rest of this section. 

\prop{6.1} Let $\phi$ be the harmonic immersion of proposition 4.4; we suppose that $\phi$ is injective. Let $j_0,\dots,j_{l-1}\in\{ 1,2,3 \}$ and let $r\not=s\in\{ A,B,C \}$. Then, 
$\phi\circ F_{j_0\dots j_{l-1}}\circ\g_{rs}$ is a curve minimising length among all curves which connect $\phi\circ F_{j_0\dots j_{l-1}}(r)$ with $\phi\circ F_{j_0\dots j_{l-1}}(s)$ and have image in the cell 
$\phi\circ F_{j_0\dots j_{l-1}}(G)$. 

\proof For simplicity, we suppose that $l=0$, i.e. that the map $F_{j_0\dots j_{l-1}}$ is the identity; the general case follows applying the argument below to the cell 
$F_{j_0\dots j_{l-1}}(G)$. 

We must prove that $\phi\circ\g_{r,s}$ is a geodesic in $\phi( G)$; to fix ideas, we suppose that $r=A$, $s=B$ and thus $\g_{r.s}=\g_{A,B}$. Since $\phi$ is injective, $\tilde A$, $\tilde B$ and 
$\tilde C$ do not lie on the same free geodesic of $H_+$: otherwise, by (4.8) $\phi(G)$ would be contained in that geodesic and then we could find, by the intermediate value theorem, $t,s\in(0,1)$ such that $\phi\circ\g_{A,B}(s)=\phi\circ\g_{A,C}(t)$, contradicting injectivity. 

Recall also that 
$\tilde A$, $\tilde B$, $\tilde C$ are oriented in the clockwise sense. Thus, applying an automorphism to $H_+$, we can also suppose that 
$$\tilde A=i,\qquad \tilde B=ai
\txt{with} a>1,\qquad{\rm Re}(\tilde C)>0  .  \eqno (6.2)$$
The idea is the following. We shall define a closed set $\bar R\subset H_+$ which is bounded on the left by $\phi\circ\g_{A,B}([0,1])$; a standard connectedness argument shows that 
$\phi(G)\subset\bar R$. Next, among the curves $\fun{\bar\eta}{[0,1]}{H_+}$ which connect $\tilde A$ with $\tilde B$ and have image in $\bar R$ we find one with minimal energy; this is a classical "geodesic with obstacle" problem; we are going to see that the infimum is finite and that it is attained on a curve $\eta$. To end the proof, we show that $\eta([0,1])=\phi\circ\g([0,1])$. Since $\eta$ is minimal among all curves which have image in $\bar R$, a larger class than curves which have image in $\phi(G)$, we get that 
$\phi\circ\g$ is a geodesic in $\phi(G)$. Now to the rigorous proof. 

\noindent{\bf Definition of the set $R$.} Formulas (6.2) and (4.8) imply that 
$\phi(G)$ does not intersect the following subset of the imaginary axis: 
$$U\colon =\{
it\st t\in(0,1)\cup(a,+\infty)
\}  .  $$
We consider the closed set of $H_+$ 
$$S\colon=\{
it\st t\in(0,1]\cup[a,+\infty)
\}\cup
\phi\circ\g_{ A, B}([0,1]) . \eqno (6.3)$$
By (6.2) we have $\phi\circ\g_{ A, B}(0)=\tilde A=i$, $\phi\circ\g_{ A, B}(1)=\tilde B=ai$; thus, $S$ is the image of a continuous curve connecting $0$ with infinity; we assert that this curve has no self-intersections. First we note an obvious fact, i.e. that $U$ has no self-intersections. Second, 
$\phi\circ\g_{A,B}([0,1])\subset\phi(G)$ does not intersect the set $U$ by (4.8). Lastly, we show that $\phi\circ\g_{A,B}$ has no self-intersections. We are supposing that $\phi$ is injective; since 
$\g_{A,B}$ is injective because it is a geodesic of $G$, we get that $\phi\circ\g_{A,B}$ is injective; these facts prove the assertion. 

As a consequence, $S$ divides $H_+\setminus S$ into two connected components, both unbounded; we call $R$ the one which contains $\tilde C+s$ for $s>0$ large; we assert that $\tilde C\in R$. Indeed, the half-line $\{ \tilde C+s\st s\ge 0 \}$ is connected and, by the last formula of (6.2) and (4.8), 
$\{ \tilde C+s\st s\ge 0 \}$ does not intersect $S$; as a consequence, this set lies in one of the two connected components of $H_+\setminus S$, which must be $R$. 

Now we note that $\phi(G)\subset\bar R$. Since $\phi$ is injective and 
$\phi(G)\subset T(\tilde A,\tilde B,\tilde C)$ by (4.8), we get that 
$\phi(G\setminus\g_{\tilde A,\tilde B}([0,1]))$ does not intersect $S$; since $\phi$ is continuous and $G\setminus\g_{\tilde A,\tilde B}([0,1])$ is connected, its image is in one of the two connected components of $H_+\setminus S$; since we just saw that $\tilde C\in R$, we get that 
$\phi(G)\subset\bar R$.

\noindent{\bf The geodesic with obstacle.} In the next lemma we show that there is a shortest curve in $H_+$ which connects $\tilde A$ with $\tilde B$ and has image in $\bar R$; before stating it, we need a definition. 

\vskip 1pc

\noindent{\bf Definitions.} We define $\rc$ as the set of all curves $\eta$ which connect $\tilde A$ with $\tilde B$ and such that $\eta([0,1])\subset\bar R$. More precisely, a curve 
$$\fun{\eta}{[0,1]}{H_+}$$is in $\rc$ if the following three points hold. 

\noindent 1) $\eta\in H^1((0,1),H_+)$; by the properties of the Sobolev space $H^1$, $\eta$ has a representative in $C([0,1],H_+)$, which justifies the following two points. 

\noindent 2) $\eta(0)=\tilde A$, $\eta(1)=\tilde B$ . 

\noindent 3) $\eta([0,1])\subset\bar R$. 

If $\eta(t)=\eta_1(t)+i\eta_2(t)\in H^1((0,1),H_+)$ we define as usual its hyperbolic energy as 
$$\int_0^1\frac{1}{\eta_2(t)^2}|\dot\eta(t)|^2\dr t  .  $$
For the obstacle problem, we define 
$$L\colon=\inf_{\eta\in\rc}
\int_0^1\frac{1}{\eta_2(t)^2}|\dot\eta(t)|^2\dr t  .  $$

\lem{6.2} The infimum above is a minimum, i.e. there is $\eta\in\rc$ such that 
$$L=\int_0^1\frac{1}{\eta_2(t)^2}|\dot\eta(t)|^2\dr t  .  \eqno (6.4)$$

\proof We begin to prove that $L$ is finite. We consider two free geodesics of $H_+$: 
$\eta_{\tilde A,\tilde C}$ is the one with $\eta_{\tilde A,\tilde C}(0)=\tilde A$, 
$\eta_{\tilde A,\tilde C}(1)=\tilde C$ while $\eta_{\tilde C,\tilde B}$ satisfies 
$\eta_{\tilde C,\tilde B}(0)=\tilde C$, $\eta_{\tilde C,\tilde B}(1)=\tilde B$. They both have finite energy, implying that the absolutely continuous curve 
$$\fun{\eta}{[0,1]}{H_+}$$
$$\eta(t)=\left\{
\eqalign{
\eta_{\tilde A,\tilde C}(2t)&\txt{if}0\le t\le\2\cr
\eta_{\tilde C,\tilde B}(2t-1)&\txt{if}\2\le t\le 1
}
\right.  $$
has finite energy. Moreover, formula (4.8) implies that $\eta\in\rc$. A consequence of these two facts is that $L$ is finite. In particular, there is a minimising sequence which, by standard results, converges to a curve $\g=\g_1+i\g_2\in\rc$ which satisfies (6.4). 

\fin

\noindent{\bf The geodesic with obstacle coincides with $\phi\circ\g_{A,B}([0,1])$.} We briefly explain the idea of the proof of lemma 6.3 below. Let $\eta$ minimise in (6.4) and let us suppose by contradiction that there is a maximal non-empty interval $(t_1,t_2)$ such that 
$\eta(t_1,t_2)\cap\phi\circ\g_{A,B}([0,1])=\emptyset$. As a consequence, $\eta(t_1,t_2)$ does not touch the obstacle and thus it is a free geodesic; up to applying an automorphism, we can suppose that $\eta(t_1,t_2)$ is a segment on the imaginary axis; since $(t_1,t_2)$ is maximal, we have that 
$\eta(t_i)=\phi\circ\g_{A,B}(s_i)$ for $i\in\{ 1,2 \}$. Let us suppose, to fix ideas, that 
$\phi\circ\g_{A,B}(s_1,s_2)$ intersects the left hand half plane.  We project onto the imaginary axis the part of $\phi(G)$ which stays to the right of $\phi\circ\g_{A,B}(s_1,s_2)$ and to the left of the imaginary axis and we show that the energy of the projection is strictly smaller that that of the minimal $\phi$, obtaining a contradiction as in lemma 4.3. 

\lem{6.3} Let $\eta\in\rc$ minimise in (6.4); then 
$$\eta([0,1])=\phi\circ\g_{A,B}([0,1])  .  \eqno (6.5)$$

\proof  We suppose by contradiction that (6.5) does not hold; since $\g$ and $\phi\circ\g_{A,B}$ are continuous, 
$$\{
t\st \eta(t)\not\in\phi\circ\g_{A,B}([0,1])
\}   \eqno (6.6)$$
is an open set; we consider one of its connected components $(t_1,t_2)\subset[0,1]$: 
$$\eta(t)\not\in\phi\circ\g_{A,B}([0,1])\txt{if}t\in(t_1,t_2)
\txt{but}\eta(t_1),\eta(t_2)\in\phi\circ\g_{A,B}([0,1])  .  \eqno (6.7)$$
Since $\eta$ is a geodesic, albeit with obstacle, it is injective and thus $\eta(t_1)\not=\eta(t_2)$. Since also $\phi\circ\g_{A,B}$ is injective, there are $s_1\not=s_2$ such that 
$\eta(t_i)=\phi\circ\g_{A,B}(s_i)$ for $i\in\{ 1,2 \}$. 

\noindent{\bf The geodesic is free.} We assert that $\eta|_{(t_1,t_2)}$ is a free geodesic in $H_+$; since $\eta\in\rc$, this is true if $\eta((t_1,t_2))$ does not intersect the set $S$ of (6.3). We know by (6.7) that $\eta((t_1,t_2))$ does not intersect $\psi\circ\g_{A,B}([0,1])$; we must prove that it does not intersect $i(0,1]\cup i[a,+\infty)$. Let us suppose by contradiction that there is $t_0\in(t_1,t_2)$ such that $\eta(t_0)\in i(0,1]\cup i[a,+\infty)$. First of all, we can exclude that $\eta(t_0)=i$ and that 
$\eta(t_0)=ia$: since these points also belong to $\phi\circ\g_{A,B}([0,1])$, in this case $(t_1,t_2)$ would not be a connected component of (6.6). Thus, $\eta(t_0)\in i\R$; let $[\a,\b]$ be the maximal interval which contains $t_0$ and such that 
$$\eta(t)\in i\R\txt{for}t\in[\a,\b]  .  $$
We have that $[\a,\b]\subset [t_1,t_2]$ and it is easy to see that the inclusion is proper unless 
$\eta([t_1,t_2])$ is contained in the imaginary axis, but in this case the assertion follows trivially because the imaginary axis is a geodesic. 

We work out the case $\b<t_2$, since the case $\a>t_1$ is analogous. Since $\eta\in\rc$ and $\b$ is maximal, for $\d>0$ small we get that 
$${\rm Re}\g(t)=0\txt{for}t\in[\a,\b]
\txt{and}{\rm Re}\eta(\b+\d)>0  .   $$
Let $\tilde\eta$ be the free hyperbolic geodesic connecting $\eta(\a-\d)$ with $\eta(\b+\d)$; the last formula easily implies that $\tilde\eta$ has image in $\bar R$. Since $\eta|_{[\a-\d,\b+\d]}$ is clearly not a geodesic of $H_+$, $\tilde\eta$ is shorter, contradicting the minimality of $\eta$ in $\rc$. 

\noindent{\bf The embedding does not minimise.} We reach the contradiction that the embedding 
$\psi$ does not minimise the energy of (5). 

Let $(t_1,t_2)$ be as in (6.7); we apply another automorphism of $H_+$ in order to move 
$\eta(t_1)$, $\eta(t_2)$ on the imaginary axis; in other words, there is $a>0$ such that, for  $s_1,s_2\in[0,1]$, 
$$\eta(t_1)=i=\phi\circ\g_{A,B}(s_1),\qquad
\eta(t_2)=ia=\phi\circ\g_{A,B}(s_2)    .    \eqno (6.8)$$
We forego to give an explicit name to the automorphism but from now on by $\tilde A$, $\tilde B$, 
$\tilde C$ and $R$ we shall mean their image by this map. 

We saw above that $\eta|_{[t_1,t_2]}$ is a free geodesic; together with the last formula this implies that $\eta([t_1,t_2])$ is contained in the imaginary axis. We suppose that $\phi\circ\g_{A,B}(s_1,s_2)$ intersects the open left hand half plane; the argument for the right hand plane is analogous. 

To fix ideas we suppose that $s_1<s_2$ and we consider the loop $\tilde\eta$ obtained joining $\eta$ with $\phi\circ\g_{A,B}$: 
$$\tilde\eta(t)=\left\{
\eqalign{
\eta(t)&\txt{if}t\in[t_1,t_2]\cr
\phi\circ\g_{A,B}(t_2-t+s_2)&\txt{if}t\in[t_2,t_2+s_2-s_1]  .   
}
\right.   $$
Formula (6.8) implies that this curve is continuous and that it is closed, i.e. 
$\tilde\eta(t_1)=\tilde\eta(t_2+s_2-s_1)$. Moreover, $\tilde\eta$ has no self-intersections. Indeed, 
$\eta|_{[t_1,t_2]}$ has no self-intersections because it is a free geodesic of $H_+$; 
$\phi\circ\g_{A,B}|_{[s_1,s_2]}$ has no self -intersections because $\phi$ is injective and 
$\g_{A,B}$ is a geodesic of $G$. Lastly,
$$\eta(t_1,t_2)\cap\phi\circ\g_{A,B}(s_1,s_2)=\emptyset$$
by the definition of $(t_1,t_2)$ in (6.7). 

By Jordan's closed curve theorem, $H_+\setminus\tilde\eta([t_1,t_2+s_2-s_1])$ has two connected components and we call $\tilde\oc$ the bounded one; part of its boundary, $\eta(t_1,t_2)$, is contained in $R$, part of it, $\phi\circ\g_{A,B}(s_1,s_2)$, is contained in $\partial R$. Since $\tilde\oc$ is a connected set which intersects $R$ only at its boundary and contains points of $R$, we get that 
$\tilde\oc\subset R$.

Since $\phi\circ\g_{A,B}(s_1,s_2)$ intersects the left hand plane we have that 
$$\oc\colon=\tilde\oc\cap\{ {\rm Re}(z)<0 \} \not=\emptyset.      $$
We set 
$$\uc=\phi^{-1}(\oc)  .  $$
We refer the reader to figure 2 below: the free geodesic $\eta|_{[t_1,t_2]}$ is the bold segment on the imaginary axis, the curve $\phi\circ\g_{A,B}$ is the curved line and the set $\oc$ is shaded. 

\includegraphics[scale=0.5]{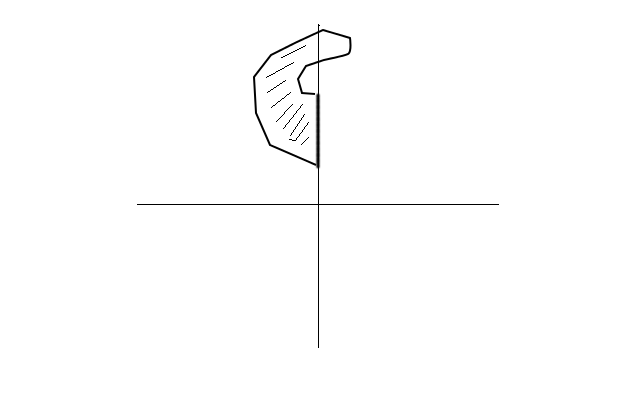}

\centerline{Figure 2}

\vskip 1pc

Now we note that $\uc$ is not empty; indeed, let us consider the points of $G$ in a thin strip around 
$\g_{A,B}(s_1,s_2)$. Since $\phi$ is continuous, these points are sent by $\phi$ into a thin strip around $\phi\circ\g_{A,B}(s_1,s_2)$; since $\tilde\oc\subset R$ it is not hard to show that this strip intersect $\oc$, implying that $\uc$ is not empty.

For the function $g$ of (4.5) we define $\tilde\psi\in\ac$ by
$$\tilde\phi(x)=\left\{
\eqalign{
\phi(x)&\txt{if}x\not\in\uc\cr
g\circ\phi(x)&\txt{if}x\in\uc  .  
}
\right.  $$

We assert that the three points below hold. 

\noindent 1) $\tilde\phi\in\dc(\dr)\oplus i\dc(\dr)$.

\noindent 2) $\tilde\phi(A)=\tilde A$, $\tilde\phi(B)=\tilde B$, $\tilde\phi(C)=\tilde C$. 

\noindent 3) $\tilde\ec(\tilde\phi)<\tilde\ec(\phi)$. 

By the first two points $\tilde\phi\in\ac$ while the third one contradicts the minimality of $\phi$, ending the proof of lemma 6.4. 

We  forego the proof of point 2), which follows easily from the definition of $\tilde\phi$ and is the same as in lemma 4.3. We also skip the proof of point 3), since again it is analogous to that of lemma 4.3; we prove point 1). 

Since the operator $\dr$ is local, it suffices to show that, for all $x_0\in G$, there is $r>0$ such that 
$\dr\tilde\phi$ exists in $B(x_0,r)$; a more formal definition would be that, for all 
$C^1$ functions $\psi$ supported in such a ball, $\tilde\phi\cdot\psi\in\dc(\dr)$. 

There are various cases. The first one is when $x_0=\phi^{-1}(z_0)$ and 
$B(z_0,r)\subset\oc^c$. In this case, for $r^\prime$ small, $\tilde\phi=\phi$ on 
$B(x_0,r^\prime)\subset\phi^{-1}(B(z_0,r))$ and we already know that $\phi$ is in the domain of $\dr$. 

The second case is when $B(z_0,r)\subset\oc$. We set $\phi=u+iv$ and we see that, for $r^\prime$ small, on $B(x_0,r^\prime)\subset\phi^{-1}(B(z_0,r))$ we have $\tilde\phi=iv$, and thus it is in the domain of $\dr$. 

The same argument applies when $z_0\in\partial\oc$ but $z_0$ is not on the imaginary axis. 

The last case is when $\phi(z_0)\in\partial\oc$ is on the imaginary axis. In this case, 
$\tilde\phi=g\circ\phi$ in $B(x_0,r^\prime)\subset\phi^{-1}(B(z_0,r))$, and the same argument of lemma 4.2 shows that $\tilde\phi$ is in the domain of $\dr$. 

\fin

\noindent{\bf Proof of point 3) of theorem 1.} We follow [16]. Let $\g$ be a geodesic in $\phi(G)$ with $\g(0)=x$ and $\g(1)=y$. It suffices to show that $\phi(G)$ is a join of curves of the type 
$\phi\circ F_{j_0\dots j_{l-1}}\circ\g_{i,j}$ with $i\not=j\in\{ A,B,C \}$. 

For a fixed $l$, we consider the cells $\phi\circ F_{j_0\dots j_{l-1}}(G)$. Since $G$ is compact and we are supposing that $\phi$ is injective, $\phi$ is a homeomorphism on its image, and thus for the interior part we have 
$$[\phi\circ F_{j_0\dots j_{l-1}}(G)]^\circ=
\phi[F_{j_0\dots j_{l-1}}(G)^\circ]  .  \eqno (6.9)$$
We define the set $\{ (t_1^\a,t_2^\a) \}_\a$ of maximal intervals such that 
$$\g((t_1^\a,t_2^\a))\subset \phi\circ F_{j_0\dots j_{l-1}}(G)^\circ    \eqno (6.10)$$
for some $j_0\dots j_{l-1}$. 

Since $(t_1^\a,t_2^\a)$ is maximal, (6.9) and (6.10) imply that 
$\g(t_1^\a),\g(t_2^\a)\in F_{j_0\dots j_{l-1}}(\{ A,B,C \})$; since $\g$ is a geodesic, it has no self-intersections and thus $\g(t_1^\a)\not=\g(t_2^\a)$. By (6.10) and the fact that $\g$ is a geodesic, 
$\g|_{[t_1^\a,t_2^\a]}$ is minimal among all curves with image in $\phi\circ F_{j_0\dots j_{l-1}}(G)$. To fix ideas, we suppose that $\g(t_1^\a)=\phi\circ F_{j_0\dots j_{l-1}}(A)$, 
$\g(t_2^\a)=\phi\circ F_{j_0\dots j_{l-1}}(B)$; by proposition 6.2 and using the notation of (6.1), for 
$t\in[t_1^\a,t_2^\a]$, up to reparametrisation we have 
$$\g(t)=\phi\circ F_{j_0\dots j_{l-1}}\circ\g_{A,B}\left(
\frac{t-t_1^\a}{t_2^\a-t_1^\a}
\right)    $$
i.e. on these intervals $\g$ is the image of a geodesic of the harmonic gasket, as we wanted. 

Recall that we fixed $l$; now we note that 
$$[0,1]\setminus\bigcup_\a(t^\a_1,t^\a_2)=[0,\e_1]\cup[1-\e_2,1]$$
with $\e_1,\e_2\tends 0$ as $l\tends+\infty$. Thus the thesis follows letting $l\tends+\infty$. 

\fin

\vskip 2pc
\centerline{\bf References} 



\noindent [1] P. Alonso-Ruiz, F. Baudoin, Dirichlet forms on metric measure spaces as Mosco limits of Korevaar-Schoen energies, preprint 2023. 

\noindent [2] P. Alonso-Ruiz, U. Freiberg, J. Kigami, Completely symmetric resistance forms on the Stretched Sierpinski Gasket, Journal of Fractal Geometry, {\bf 5-3}, 227-277, 2018. 

\noindent [3] A. Ambrosetti, A. Malchiodi, Nonlinear Analysis and semilinear elliptic problems, Cambridge, 2006.





\noindent [4] L. Ambrosio, N. Gigli, G. Savar\'e, Heat flow and calculus on metric measure spaces with Ricci curvature bounded below - the compact case. Analysis and numerics of Partial Differential Equations, 63-115, Springer, Milano, 2013.  

\noindent [5] F. Baudoin, Korevaar-Schoen-Sobolev spaces and critical exponents in metric measure spaces, preprint 2024.



\noindent [6] U.Bessi, Another point of view on Kusuoka's measure, DCDS, {\bf 41-7}, 3251-3271, 2021.  


\noindent [7] U. Bessi, Harmonic embeddings of the stretched Sierpinski gasket, NoDEA, {\bf 30-6}, paper no. 80, 34pp.  

\noindent [8] U. Bessi, Hodge theory on the harmonic gasket and other fractals, to appear o Nolinear Analysis, TMA. 

\noindent [9] G. Bouchitt\'e, G. Buttazzo, P. Seppecher, Energies with respect to a measure and applications to low dimensional structures, Calc. Var. and PDE, {\bf 5}, 37-54, 1997. 


\noindent [10] S. Cao, H. Qiu, A topological proof of the non-degeneracy of harmonic structures on Sierpinski gaskets, Analysis in theory and applications, {\bf 36}, 510-516, 2020.

\noindent [11] J. Eels, J. H. Sampson, Harmonic mappings of Riemannian manifolds, American Journal of Mathematics, {\bf 86}, 109-160, 1964. 

\noindent [12] A. Johansson, A. \"Oberg, M. Pollicott, Ergodic theory of Kusuoka's measures, J. Fractal Geom., {\bf 4}, 185-214, 2017.  

\noindent [13] J. Heinonen, P. Koskela, N. Shanmugalingam, J. T. Tyson, Sobolev spaces on metric measure spaces, Cambridge, 2015. 

\noindent [14] N. Kajino, Analysis and geometry of the measurable Riemannian structure on the Sierpiski gasket, Contemporary Mathematics, {\bf 600}, 91-133, 2013. 


\noindent [15] J. Kigami, Analysis on fractals, Cambridge, 2001. 

\noindent [16] J. Kigami, Measurable Riemannian geometry on the Sierpinski gasket: the Kusuoka measure and the Gaussian heat kernel estimates, Math. Ann., {\bf 340}, 781-804, 2008.


\noindent [17] P. Koskela, Y. Zhou, Geometry and Analysis of Dirichlet forms, Adv. Math., 
{\bf 231}, 2755-2801, 2012.



\noindent [18] S. Kusuoka, Dirichlet forms on fractals and products of random matrices, Publ. Res. Inst. Math. Sci., {\bf 25}, 659-680, 1989. 






\noindent [19] I. D. Morris, Ergodic properties of matrix equilibrium state, Ergodic Theory and Dyn. Sys., {38/6}, 2295-2320, 2018.






\noindent [20] M. Piraino, The weak Bernoulli property for matrix equilibrium states, Ergodic theory Dynam. Systems, {\bf 40}, 2219-2238, 2020.







\noindent [21] R. S. Strichartz, Differential equations on fractals, a tutorial, Princeton, 2006. 

\noindent [22] M. E. Taylor, Partial differential equations, Basic theory, Springer, Berlin, 1996

\noindent [23] A. Teplyaev, Harmonic coordinates on fractals with finitely ramified cell structure, Canad. J. Math., {\bf 60}, 457-480, 2008.

\noindent [24] A. Teplyaev, Energy and Laplacian on the Sierpinski gasket, Proceedings of symposia in pure Mathematics, {\bf 72-1}, 131-154, Providence, R. I., 2004. 

\noindent [25] K. Tsougas, Non-degeneracy of the harmonic structure on Sierpinski gaskets, J. Fractal Geom., {\bf 6}, 143-156, 2019. 




\end